\documentclass{amsart}
\usepackage[english]{babel}
\usepackage{graphicx}
\usepackage{amssymb}
\usepackage{amsmath}
\usepackage[foot]{amsaddr}
\usepackage{amsfonts}
\usepackage{bbm}
\usepackage{xcolor}
\usepackage{float}
\usepackage{algorithm}
\usepackage{algpseudocode}
\usepackage{rotating}
\usepackage[normalem]{ulem}

\usepackage{amsmath}
\usepackage{amssymb}
\usepackage{mathrsfs}
\usepackage{color}
\usepackage{graphicx}
\usepackage{overpic}
\usepackage{url}
\usepackage{subcaption}
\usepackage{enumerate}
\usepackage{amsthm}
\usepackage{moreverb}
\usepackage[pdftex,colorlinks,bookmarksopen,bookmarksnumbered,citecolor=red,urlcolor=red]{hyperref}

\def\0{\boldsymbol{0}}

\usepackage{tikz}
\usepackage[top=1in, bottom=1.25in, left=1.25in, right=1.25in]{geometry}

\usepackage{xcolor}
\usepackage{hyperref}

\newcommand{\bE}{\mathbf E}

\newcommand{\bQ}{\mathbf Q}
\newcommand{\bR}{\mathbf R}

\newcommand{\bx}{\mathbf x}

\begin{document}

\title[IREE algorithm for mesh adaptation with FV]{{An isotropic recovery-based error estimator algorithm for mesh adaptation in a finite volume environment with application to atmospheric flows}} 

\author{L. Besabe$^1$, M.~ Girfoglio$^2$, S.~Perotto$^3$,  A.~Quaini$^4$, G.~Rozza$^5$}

\address{$^1$School of Mathematical and Statistical Sciences, Clemson University, Clemson, SC 29634-0975, USA; lbesabe@clemson.edu}

\address{$^2$ Department of Engineering, University of Palermo, Viale dlle Scienze Ed. 7, 90128, Palermo, Italy; michele.girfoglio@unipa.it}

\address{$^3$ MOX-Department of Mathematics – Politecnico di Milano, 
Piazza Leonardo da Vinci 32, I-20133 Milano, Italy; simona.perotto@polimi.it
}

\address{$^4$ Department of Mathematics, University of Houston, 3551 Cullen Blvd, Houston TX 77204, USA; aquaini@central.uh.edu}

\address{$^5$ mathLab, Mathematics Area, SISSA, via Bonomea 265, I-34136 Trieste, Italy; grozza@sissa.it 
}

\begin{abstract}
We develop an Isotropic Recovery-based Error Estimator (IREE) to drive mesh adaptation within a finite volume framework. Recovery-based error estimators are widely used in practice thanks to their simplicity, relying solely on the available discrete solution and a suitable post-processing step. While recovery-based estimators are well established in the finite element framework, their application to finite volume methods, despite the widespread use of the latter in industrial and commercial codes, remains largely unexplored and motivates the present study.
We assess the performance of the proposed IREE-driven mesh adaptation procedure through well-known benchmarks for dry atmospheric flows modeling and compare it  against a widely-used plain mesh adaptation algorithm. Both qualitative and quantitative results show that 
the IREE algorithm is more effective at suppressing numerical instabilities compared with the plain approach and, for most of the simulation time, provides a superior accuracy, though at a moderately higher computational cost. In addition, both the adaptive approaches offer substantial computational savings (between $35$\% and $94$\%) relative to simulations based on fixed uniformly fine meshes.
These findings demonstrate that IREE-based mesh adaptation is a promising and effective strategy for atmospheric flow simulations and has the potential to substantially reduce the computational effort of generating reanalysis-quality data.
\end{abstract}

\maketitle

\section{Introduction}\label{sec:intro}
Atmospheric flow simulations remain among the most computationally demanding
applications in scientific computing. Accurate weather prediction and climate
modeling require the numerical resolution of multiscale phenomena spanning a
wide range of spatial and temporal scales, from planetary circulation patterns
to localized structures associated with fronts, gravity waves, convection, and
complex topography. Despite continuous advances in numerical methods and
high-performance computing, the generation of high-fidelity simulation and
reanalysis datasets remains computationally expensive. 
This computational burden has motivated several research directions. In
particular, recent machine-learning-based forecasting systems
\cite{pathak2018model,rasp2021data,schultz2021can,weyn2019can,GraphCast,FourCastNet,Pangu-Weather}
have demonstrated remarkable predictive capabilities. Most of these approaches,
however, rely heavily on large quantities of high-quality training data,
typically obtained from expensive numerical simulations and reanalysis products.
Consequently, reducing the computational cost required to generate such datasets
remains an important challenge.

Mesh adaptation offers a natural strategy to address this issue. By dynamically
concentrating computational resources only where they are needed, adaptive
meshes can significantly improve the accuracy-to-cost ratio of numerical
simulations. While mesh adaptation is nowadays routinely employed in several areas of computational fluid dynamics and is available in many industrial and commercial software packages, its penetration into mainstream atmospheric modeling remains comparatively limited. This is somewhat surprising, since atmospheric flows exhibit exactly the type of localized multiscale structures for which adaptive techniques are expected to be particularly effective. 
Nevertheless, mesh adaptation for atmospheric applications has attracted sustained research interest over the last several decades. Despite a substantial body of literature and a number of successful implementations, adaptive strategies have not yet achieved the same level of adoption in atmospheric modeling as in many other CFD applications. Early contributions focused on adaptive strategies for
hyperbolic conservation laws and numerical weather prediction
\cite{BERGER1984484,SKAMAROCK198927,Skamarock1993}, while the first operational
atmospheric model employing horizontally adaptive grids was introduced in
\cite{Bacon2000}. Subsequent developments combined adaptive meshes with
high-order Galerkin discretizations
\cite{Chen2011,MULLER2013371,KOPERA201492,KOPERA201590}
and demonstrated the potential of mesh adaptation for the simulation of
atmospheric flows over complex topography
\cite{Li2021,YAMAZAKI2022111217,Orlando2024}
and highly localized events such as tropical cyclones
\cite{Hendricks2016,Tissaoui2025}.

A common feature of most existing atmospheric adaptation strategies is that
adaptation is performed predominantly in the horizontal directions, while the
vertical mesh remains structured. This choice is largely motivated by the need
to preserve compatibility with column-based physical parameterizations used to
model clouds and precipitation. More recently, however, vertically adaptive
approaches \cite{Carstensen2025} and fully unstructured adaptive meshes
\cite{Tissaoui2023} have demonstrated that these limitations can be relaxed
without compromising the accuracy of atmospheric simulations.

The criteria driving mesh adaptation are equally diverse. Existing approaches
include Richardson-based estimates of the local truncation error
\cite{BERGER1984484}, adaptation driven by selected solution variables or by
their gradients
\cite{MULLER2013371,KOPERA201492,KOPERA201590,Hendricks2016,Chen2021,Li2021,Tissaoui2025},
optimal-transport-based relocation strategies
\cite{YAMAZAKI2022111217}, and, more recently, data-driven approaches relying
on machine learning \cite{Gan2026}. In several of the references above
\cite{MULLER2013371,KOPERA201492,KOPERA201590,Hendricks2016,Chen2021,Li2021,Tissaoui2025},
mesh refinement and coarsening are driven directly by selected solution fields
or by norms of their gradients. Owing to its widespread adoption, we refer to
this family of approaches as Plain Mesh Adaptation (PMA) and use it as
the reference methodology for comparison throughout this work.

Despite this substantial body of work, comparatively little attention has been
devoted to rigorous a posteriori error estimation techniques in the
context of atmospheric mesh adaptation. A posteriori approaches exploit the
computed numerical solution to identify regions where the discretization error
is larger and additional resolution is required. Their main strength lies in the
availability of computable error estimates, which can be used to drive the
adaptive process in a mathematically informed manner. The adaptive procedure is
therefore inherently iterative, involving a continuous feedback loop between the
numerical solution, the error estimate, and the computational mesh.

Among the many available a posteriori error estimation strategies,
recovery-based error estimators are particularly attractive due to their
simplicity and their ability to estimate the discretization error using only the
available discrete solution. The key idea is to replace the unavailable exact
quantity appearing in the error definition with a suitably reconstructed
(recovered) approximation computed from the discrete solution itself \cite{Zienkiewicz1987,ZZ92_parte1,ZZ92_parte2}.
Recovery-based error estimators have been
extensively investigated within the finite element framework and successfully
applied to a broad range of problems
\cite{Markowich1988,Mu2013,Yan2001}. The anisotropic counterpart of recovery-based error estimators was introduced in \cite{MichelettiPerotto2010} to combine the simplicity and efficiency of recovery procedures with the superior resolution capabilities of anisotropic mesh adaptation. The resulting methodology has since been successfully applied to a large variety of computational problems \cite{Porta2012,Esfandiar2015,Artina15,Soli2019,
Chiappa19,FMP20,Belz21,reset1,Mauri24,Temellini25}.

Despite the generality of recovery-based estimators, the vast majority of the
literature develops and applies them within finite element settings. Their
extension to finite volume discretizations, which remain the method of choice in
many industrial and commercial CFD codes, has received considerably less
attention. At the same time, recovery-based adaptive strategies have seen only
limited use in atmospheric flow simulations. The combination of these two
aspects constitutes the main motivation for the present work. Within this setting, we deliberately restrict our attention to isotropic mesh
adaptation. While anisotropic strategies can provide additional gains in the
presence of strongly directional solution features, they introduce further
algorithmic and mesh-management complexities. Given the novelty of combining
recovery-based error estimation, finite volume discretizations, and atmospheric
flow simulations, we regard the isotropic framework as the most appropriate
starting point for assessing the effectiveness of the proposed methodology.

The purpose of this paper is therefore twofold. First, we develop an isotropic
recovery-based error estimator (IREE) algorithm for mesh adaptation in a finite
volume environment. Second, we assess its effectiveness through a collection of
well-established atmospheric benchmarks. The proposed adaptive strategy is
systematically compared against PMA and shown to substantially reduce the
computational cost of atmospheric simulations while preserving the desired level
of accuracy. In particular, the proposed methodology achieves reductions of up
to 91\% in computational time, highlighting its potential for the efficient
generation of high-quality atmospheric datasets.

The remainder of the paper is organized as follows.
Section~\ref{sec:IREE} describes the proposed IREE mesh adaptation algorithm and
introduces PMA for comparison purposes.
Section~\ref{sec:pbd} presents the compressible Euler equations for low-Mach
stratified flows together with the adopted spatial and temporal
discretizations.
Numerical results are reported in Sec.~\ref{sec:flatBench}, while concluding
remarks are drawn in Sec.~\ref{sec:conc}.

\section{An isotropic recovery-based mesh adaptation algorithm}\label{sec:IREE}

This section presents the proposed mesh adaptation technique based on an isotropic recovery-based error estimator (IREE)  for a generic PDE problem 
in primary variables space $\bx \in \Omega$, where $\Omega$ is a computational domain of interest, and  time $t \in (0, T_f]$. For the time discretization of such PDE problem, we introduce a time step $\Delta t = T_f/N_t$ for given $N_t\in \mathbb N$, and set $t^n =n\Delta t $ for $n = 0, \dots, N_t$, with $t^0=0$. We
denote with $y^n$ the approximation of a scalar quantity 
$y$ at time $t^n$. For space discretization,
we adopt a finite volume method. For this, let us consider a Cartesian partition, ${\mathcal T}_h$, of the computational domain $\Omega$ into cells or control volumes. We will denote with $K$ the generic cell in the mesh, with $N_c$ the total number of cells, 
and with $N_V$ the total number of vertices.

Let $\bQ(t)$ be the gradient of a  
time dependent solution field of interest. This quantity could be any physical field computed in the simulation and so it could be a scalar or a vector. 
We drive the mesh adaptation
by an a posteriori isotropic recovery-based error estimator for the $L^2$-norm of  error $\bE_h^n= \bQ_h^n - \bQ(t^n)$, where $\bQ_h^n$ is the approximation of the exact solution $\bQ(t^n)$, at time $t^n$, given by a finite volume method. 
Obviously, we cannot quantify the error itself because we do not have  $\bQ(t^n)$.
Thus, we need to find a computable quantity able to estimate the error. 
Following \cite{Zienkiewicz1987,ZZ92_parte2}, the squared $L^2$-norm of $\bE_h^n$ is estimated by
\begin{equation}
    \|\bE_h^n\|^2_{L^2(\Omega)} \approx \eta^2(t^n) = \sum_{K\in\mathcal{T}_h^n} \eta_K^2(t^n) = \sum_{K\in\mathcal{T}_h^n} \| \bQ_h^n - \bR(\bQ_h^n) \|^2_{L^2(K)},
\end{equation}
where $\mathcal{T}_h^n$ is the mesh at time $t^n$, $\eta_K$ is the local contribution to the (global) error estimator $\eta$, and $\bR(Q_h^n)$ is the so-called recovered gradient, namely an approximation of $\bQ(t^n)$ obtained by using the computed $\bQ_h^n$. 

\begin{figure}[h]
    \centering
    \begin{tikzpicture}[scale=0.9]
        \draw[thick] (-1,-1) rectangle (1,1);

        \draw[thick] (-1,1) rectangle (1,3);

        \draw[thick] (-1,-3) rectangle (1,-1);

        \draw[thick] (-3,-1) rectangle (-1,1);

        \draw[thick] (1,-1) rectangle (3,1);
        \draw[thick] (-3,1) rectangle (-1,3);
        \draw[thick] (-3,-3) rectangle (-1,-1);
        \draw[thick] (1,1) rectangle (3,3);
        \draw[thick] (1,-3) rectangle (3,-1);

        \draw[thick] (2,-1) -- (2,1);   
        \draw[thick] (1,0) -- (3,0);     
        \draw[thick] (0,1) -- (0,3);
        \draw[thick] (-1,2) -- (1,2); 
        \draw[thick] (1,2) -- (3,2); 
        \draw[thick] (2,3) -- (2,1); 
        \draw[thick] (2,2.5) -- (3,2.5); 
        \draw[thick] (2.5,3) -- (2.5,2); 

        \filldraw[fill=lightgray, draw=black, thick] (-1,-1) rectangle (1,1);
        \filldraw[fill=lightgray, draw=black, thick] (-3,-1) rectangle (-1,1);
        \filldraw[fill=lightgray, draw=black, thick] (-1,-3) rectangle (1,-1);
        \filldraw[fill=lightgray, draw=black, thick] (1,0) rectangle (2,1);
        \filldraw[fill=lightgray, draw=black, thick] (1,-1) rectangle (2,0);
        \filldraw[fill=lightgray, draw=black, thick] (-1,1) rectangle (0,2);
        \filldraw[fill=lightgray, draw=black, thick] (0,1) rectangle (1,2);

        \node at (0,0) {$\bullet$};
        \node at (1.5,0.5) {$\bullet$};
        \node at (1.5,-0.5) {$\bullet$};
        \node at (-2,0) {$\bullet$};
        \node at (-0.5,1.5) {$\bullet$};
        \node at (0.5,1.5) {$\bullet$};
        \node at (0,-2) {$\bullet$};

        \node at (0.3,0.3) {$P_K$};
        \node at (-0.7,-0.7) {$K$};

    \end{tikzpicture}
    \caption{Generic cell $K$ with centroid $P_K$ and associated patch $\Delta_K$ shaded in gray. The black dots mark the cell centroids.}
    \label{fig:quadrati}
\end{figure}
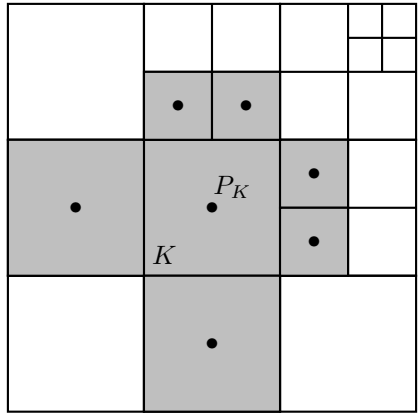

As underlined in Sec.~\ref{sec:intro}, the vast majority of the recovery-based estimators, hence the recovery operators available in the literature, are designed for finite elements methods. 
We introduce here the operator we adopt for a  finite volume approximation. To define it, let 
$\Delta_K$ be a patch of control volumes, namely,
cell $K$ itself and its neighboring cells (see Fig.~\ref{fig:quadrati} for a sketch). Then, the value of 
$\bR(\bQ_h^n)$ at the centroid of the  cell
$K$, denoted with $P_K$, is given by 
\begin{equation}\label{eq:R}
    \bR(\bQ_h^n) \big|_{P_K} = |\Delta_K|^{-1} \sum_{T\in\Delta_K}\frac{|T|}{N+1} \bQ_h^n \big|_{P_T}
\end{equation}
where $|\Delta_K|$ is the volume of the patch, 
$N+1$ is the number of cells in the patch, 
$|T|$ is the volume of cell $T$, and $\bQ_h^n\big|_{P_T}$ the value of $\bQ_h^n$ at the centroid $P_T$ of the cell $T$. This recovery operator takes inspiration from the area-weighted average formula presented in \cite{Rodri94} and merges it with ideas from
\cite{Africa2023}. 
Ultimately, \eqref{eq:R} produces an enriched
approximation of $\bQ(t^n)$ that takes into account both the numerical solution $\bQ_h^n$ and the features of the current mesh. 

The goal is to employ the IREE $\eta$ 
to generate a computational mesh able to follow the features of the solution at hand. 
To adapt an initially given mesh, we adopt
quadtree-based decomposition, a well-known hierarchical method for two-dimensional spatial domains that splits a rectangular cell (parent) into four equal rectangles (children) by halving its edges. This concept naturally extends to three dimensions with octrees, where each rectangular prism (parent) is subdivided into eight smaller rectangular prisms (children) (ee, e.g., \cite{Morton1966,Samet1985}). In particular, to start the mesh adaptation process, we need:
\begin{itemize}
\item[-] the initial mesh $\mathcal{T}_h^0$; 
\item[-] a refine interval $I_r>1$ to set every how many $\Delta t$ the mesh is adapted during the simulation;    
\item[-] the maximum level of refinement, $l_\text{max}$, allowed for each cell;
\item[-] the maximum number, $N_\text{max}$, of cells allowed in the computational domain;
\item[-] the constant parameters, $\delta_1$ and $\delta_2$, tuning the refining and coarsening procedures;
\item[-] a tolerance $tol$ on the estimated error.
\end{itemize}
Notice that, since re–meshing at every time step could be computationally too expensive,
the mesh adaptation takes place 
only if the time instant $t^n$ is a multiple of $I_r$. 

\vspace*{0.1cm}

The IREE algorithm follows a standard solve-estimate-mark-adapt paradigm as in~\cite{Africa2023}, and consists of the three main steps:
\begin{enumerate}
    \item \textbf{Initialization phase}: all the cells $K$ of the initial mesh $\mathcal{T}_h^0$ are labeled with a refinement level $l_K=0$. In general, parameter $l_K$ records the number of times cell $K$ has been split with respect to the original mesh.
    At the end of the initialization phase, the simulation is started.
    \item \textbf{Marking phase}: if $t^n$, $n \geq 1$, is a multiple of the refine interval $I_r$, then the marking process starts, otherwise, the mesh remains unchanged and the simulation continues. During the marking phase, one computes the local error estimate, $\eta_K$, for each cell. If 
    \begin{equation}\label{refineC}
    \eta_K \geq \delta_1 \frac{tol}{\sqrt{N_{c}}},
    \end{equation}
    and the refinement level does not exceed $l_\text{max}$, cell $K$ is marked
    for refinement. If
    \begin{equation}\label{coarseC}
    \eta_K\leq \delta_2\frac{tol}{\sqrt{N_{c}}},
    \end{equation} 
    cell $K$ is marked for coarsening,
    provided that its refinement level is above 0
    (i.e., we do not allow cells to be coarsened beyond the initial mesh size). If neither of the previous conditions are met, the cell remains unchanged. Let us stress that since $N_c$ varies during the simulation, the thresholds for refinement/coarsening vary with it. Moreover, criteria \eqref{refineC} and \eqref{coarseC} are intended to achieve error equidistribution over the adapted mesh. At the end of the marking phase, we have two lists of cells: one with the candidate cells for refining and the other with the candidate cells for coarsening.
    \item \textbf{Adaptation phase}: once all the cells are marked, the mesh is updated considering other possible constraints set by the user. For instance, before the refinement of the marked cells, it is verified that the predicted total number of cells would not exceed $N_\text{max}$. If the number of cells is projected to exceed $N_\text{max}$, 
    only a subset of cells gets refined. The selection criterion is based on the value of the local error estimator: the larger $\eta_K$, the higher the priority in the list.
    The selected cells are split according to a  quadtree- or octree-based decomposition.
    In order to avoid sharp changes in the mesh grading, a layer of cells, called buffer layer, is employed to “bridge” the gap between fine and coarse regions of the mesh. This buffer layer is such that the difference between the refinement levels of two adjacent cells is one. 
    If coarsening is required, a group of adjacent cells are merged together, provided that all of them are marked for coarsening.\\ 
    The solution fields in the PDE problem are then mapped onto the new mesh and, if fluxes need to be computed, they are approximated on the newly created faces.
\end{enumerate}

\vspace*{0.1cm}

Algorithm \ref{algo2} summarizes the main steps of the IREE-driven mesh adaptation strategy.

\begin{algorithm}
\caption{IREE-driven mesh adaptation algorithm}
\label{algo2}
\begin{algorithmic}[1]
    \setlength{\itemsep}{1.2mm} 
    \State $\boldsymbol{\rm input}: \mathcal{T}_h^0,\,T_{f}\,,\Delta t\,,I_r\,,\delta_1\,,\delta_2\,,tol\,,l_\text{max}$ 
    \While{$t^n\leq T_f$}
        \State compute $\bQ_h^n$ 
        \If{\texttt{mod}($t^n$,\,$I_r$)$==0$}
        \For{$K=1,2,...,N_{c}$} \Comment{Marking phase}
             \State compute $\eta_K$
             \If{$\eta_K \geq \delta_1 \frac{tol}{\sqrt{N_{c}}}\,\,$ \texttt{and} $\,\,l_K<l_\text{max}$}
                 \State mark $K$ for refinement
             \ElsIf{$\eta_K\leq \delta_2\frac{tol}{\sqrt{N_{c}}}\,\,$ \texttt{and} $\,\,l_K>0$}
                \State mark $K$ for coarsening
              \Else
                  \State $K$ remains unchanged
              \EndIf
        \EndFor
        \State update of the mesh cells and of $l_K$ \Comment{Adaptation phase}
        \EndIf
        \State $t^n = t^n +\Delta t$
    \EndWhile
\end{algorithmic}
\end{algorithm}

\subsection{A plain mesh adaptation algorithm}\label{sec:plain}
This section introduces a mesh adaptation algorithm widely used in the literature \cite{MULLER2013371,KOPERA201492,Li2021,Hendricks2016,Tissaoui2025} and in software packages (e.g., OpenFOAM\textsuperscript{\textregistered}) because of its ease of implementation and effectiveness. For simplicity, we refer to this approach as Plain Mesh Adaptation (PMA) algorithm.
Assuming that the variable driving the mesh adaptation is $\bQ_h^n$ as in the previous section, this algorithm relies on the scalar field
\begin{equation}\label{eq:alpha_k}
    \alpha_K (t^n)= \frac{\|\bQ_h^n\|_{L^2(K)}}{\max\limits_{K \in\mathcal{T}_h^n}\|\bQ_h^n\|_{L^2(K)}}
\end{equation}
which is computed for every cell $K$.
Most of the inputs required by the PMA algorithm have already been introduced in the previous section in the context of the IREE-based mesh adaptation strategy. Additionally, we need:
\begin{itemize}
    \item[-] a refinement range $\left[\alpha_\text{min}^\text{ref},\alpha_\text{max}^\text{ref}\right]\subset \left(0,1\right)$; 
    \item[-] a coarsening threshold $l_c$.
\end{itemize}
The PMA procedure has three main phases too. The first phase is the same as in the IREE algorithm. The other two are as follows:
\begin{enumerate}
    \item[(2)] \textbf{Refining phase}: 
    If the value of $\alpha_K(t^n)$ in \eqref{eq:alpha_k}
    falls within the range $\left[\alpha_\text{min}^\text{ref},\alpha_\text{max}^\text{ref}\right]$, then cell $K$ is marked for refinement provided that $l_K< l_\text{max}$. 
    Once the marking process is completed, the total number of cells in the adapted mesh is estimated, assuming that all marked cells are refined. If this number is less than $N_\text{max}$, all the marked cells are refined. Otherwise, the marked cells are ranked according to the indicator $e_K= \min\{\alpha_K(t^n) - \alpha_\text{min}^\text{ref}, \alpha_\text{max}^\text{ref} - \alpha_K(t^n)\}$
and only the first $N_\text{max}$ cells, corresponding to the largest values of $e_K$, are selected for refinement.
The selected cells are split through quadtree- or octree-based decomposition. As in the IREE-based mesh adaptation strategy, a buffer layer is introduced around the marked cells. 
    The solutions to the PDE problem are then mapped onto the new mesh and, if fluxes need to be computed, they are approximated on the newly created faces. 
    \item[(3)] \textbf{Coarsening phase}:
    Rather than checking the value of $\alpha_K$ in each cells, in this phase
    one checks the maximum of $\alpha_K$ for all the cells that share a given
    vertex $V_j$, $j = 1, \dots, N_V$. For this purpose, we call 
    $\mathcal{N}(V_j)$ the patch of cells that share $V_j$.
    If $\displaystyle \max_{K \in \mathcal{N}(V_j)} \alpha_{K}$ is smaller than $l_c$, then $V_j$ is a candidate for removal. Vertex $V_j$ is actually removed if the cell in
$\mathcal N(V_j)$ can be merged and the resulting cell does not become
coarser than the original mesh (i.e., the corresponding refinement level $l_K$
remains non-negative), while preserving the buffer layer.
\end{enumerate}
The main steps of the PMA strategy are summarized in Algorithm~\ref{algo1}.


\begin{algorithm}
\caption{PMA algorithm}
\label{algo1}
\label{algo_openfoam_unrefinement}
\begin{algorithmic}[1]
    \setlength{\itemsep}{1.3mm} 
    \State $\boldsymbol{input}: 
    \mathcal{T}_h^0,\,T_{f},\,\Delta t,\,I_r,\,l_\text{max},\,l_c,\,\alpha_\text{min}^\text{ref},\,\alpha_\text{max}^\text{ref}$    

    \While{$t^n \leq T_f$}
        \State compute ${\bf Q}_h^n$
        
        \If{\texttt{mod}($t^n$, $I_r$)$==0$}
            \For{$i=1,2,\dots,N_c$}\Comment{Refining phase}
                \State compute $\alpha_{K_i}$ in \eqref{eq:alpha_k}
                \If{$\alpha_\text{min}^\text{ref} \leq \alpha_{K_i} \leq \alpha_\text{max}^\text{ref}$ \texttt{and} $l_{K_i} < l_\text{max}$}
                    \State mark cell $K_i$ for refinement
                \EndIf
            \EndFor
            \If{the predicted number of cells is less than $N_{\max}$}
            \State refine all cells  selected for refinement
            \Else
            \State rank the marked cells according to $e_{K_i}$ and refine the first $N_{\max}$ cells
            \EndIf
            \For{$j = 1, \dots, N_V$}\Comment{Coarsening phase}
                \If{$\displaystyle \max_{K \in \mathcal{N}(V_j)} \alpha_{K} < l_c$}
                    \State mark vertex $V_j$ for removal
                \EndIf
            \EndFor
            \If{the refinement level of the cell yielded by the removal of $V_j$ is non-negative}
            \State remove $V_j$
            \EndIf
        \EndIf

        \State $t^n = t^n + \Delta t$
    \EndWhile
\end{algorithmic}
\end{algorithm}

\section{Atmospheric flow simulations: model and discretization}
\label{sec:pbd}
We consider the dynamics of dry atmosphere in a fixed spatial domain $\Omega$ by neglecting the effects of moisture, solar radiation, and heat flux from the ground. 
An established model for such dynamics is provided by the weakly compressible Euler equations. To state this model, we introduce some notation. Let
$\rho$ be the air density, $\mathbf{u} = (u, v, w)^T$ the wind velocity, and $e$ the total energy density. Note that $e = c_v T + |\mathbf{u}|^2 / 2 + gz$, where $c_v$ is the specific heat capacity at constant volume, $T$ is the absolute temperature, $g$ is the gravitational constant, $z$ is the vertical coordinate and $|\mathbf{u}|$ denotes the velocity magnitude.

Conservation of mass, momentum, and energy in terms of $\rho$, $\mathbf{u}$, and $e$ in 
domain $\Omega$ over a time interval of interest $(0, T_f]$ are given by:
\begin{align}
    &\frac{\partial \rho}{\partial t} + \nabla \cdot (\rho \mathbf{u}) = 0 \label{eq:1}\\[2mm]
    &\frac{\partial (\rho \mathbf{u})}{\partial t} + \nabla \cdot (\rho \mathbf{u} \otimes \mathbf{u}) + \nabla p + \rho g \hat{\mathbf{k}} = \mathbf{0}\label{eq:2} \\[2mm]
    &\frac{\partial (\rho e)}{\partial t} + \nabla \cdot (\rho e \mathbf{u}) + \nabla \cdot (\mathbf{u} p) = 0,\label{eq:3}
\end{align}
where $\hat{\mathbf{k}}$ is the unit vector aligned with the vertical axis $z$ and $p$ is pressure. 
To close system \eqref{eq:1}-\eqref{eq:3}, we need a thermodynamic equation of state for $p$. Assuming that dry air behaves like an ideal gas, we have:
\begin{equation}
    p = \rho R T \,,\label{eq:4}
\end{equation}
where $R$ is the specific gas constant of dry air. Moreover, for numerical stability, we add an artificial diffusion term to equations~\eqref{eq:2}-\eqref{eq:3} which become 
\begin{align}
    &\frac{\partial (\rho \mathbf{u})}{\partial t} + \nabla \cdot (\rho \mathbf{u} \otimes \mathbf{u}) + \nabla p + \rho g \hat{\mathbf{k}} - \mu_a \Delta \mathbf{u} = \mathbf{0} \label{eq:5}\\[2mm]
    &\frac{\partial (\rho e)}{\partial t} + \nabla \cdot (\rho e \mathbf{u}) + \nabla \cdot (\mathbf{u} p) - c_p \frac{\mu_a}{Pr} \Delta T = 0, \label{eq:6}
\end{align}
where $\mu_a$ is a constant (artificial) diffusivity coefficient, $Pr$ is the Prandtl number and $c_p$ is the specific heat capacity at constant pressure. This is equivalent to using a basic eddy viscosity model in Large Eddy Simulation (LES). More sophisticated LES models can be found, e.g., in \cite{Bou-Zeid2005,Marras2015,GQR_OF_clima,CGQR}.

Starting from problem (\ref{eq:1}), (\ref{eq:4})--(\ref{eq:6}), we derive the formulation that we will use for the numerical experiments. Let $p_0$ and $\rho_0$ be the pressure and density in  hydrostatic balance. 
From \eqref{eq:4} we have: 
\begin{align}
\rho_0 = \dfrac{p_0}{R T_0},
\label{eq:9}
\end{align} 
and to compute $p_0$ we solve the following hydrostatic equation:
\begin{equation}
    \nabla p_0 + \rho_0g\boldsymbol{\hat{k}} = \boldsymbol{0}. \label{eq:14}
\end{equation}
We write the pressure $p$
and density $\rho$ as the sum of fluctuations ($p'$ and $\rho'$) with respect to the values in  hydrostatic balance:
\begin{align}
p &= p_0 + p', \label{eq:p_splitCE2} \\
\rho &= \rho_0 + \rho'. \label{eq:rhoPrime} 
\end{align}
These splittings allow us to avoid numerical issues associated with the interaction of large gradients of the buoyancy force with mesh non-orthogonality. 
By plugging \eqref{eq:p_splitCE2} and \eqref{eq:rhoPrime} 
into (\ref{eq:5}), we get:
\begin{equation}
    \frac{\partial (\rho \mathbf{u})}{\partial t} + \nabla \cdot (\rho \mathbf{u} \otimes \mathbf{u}) + \nabla p' + \rho' g \boldsymbol{\hat{k}} - \mu_a \Delta \mathbf{u} = \mathbf{0}. \label{eq:13}
\end{equation}
For the numerical solution of the problem introduced above, we will employ a pressure-based method. 
This class of approaches was originally developed for weakly compressible flows. The term pressure-based refers to the fact that the pressure field is obtained by solving a pressure (or pressure-correction) equation derived from the continuity and momentum equations. The alternative is given by density-based methods (see, e.g., \cite{KAPPELI2014199,Clinco2024}), which use the continuity equation to get the density field, while the pressure field is computed from the equation of state.

In \cite{GIRFOGLIO2025106510} it was shown that, 
when using a pressure-based approach and a space discretization based on finite volumes, the formulation of the weakly compressible Euler equations is not the most accurate, reliable, and robust model for the simulation of mesoscale atmospheric flows. The modification reported below
results in an improved formulation of the problem from the numerical point of view.

Let us introduce the potential temperature $\theta$, i.e., the temperature a parcel of dry air would have if it were expanded or compressed adiabatically to the atmospheric pressure at the ground, $p_g = 10^5$ Pa, and the so-called Exner pressure, $\pi$, given by
\begin{equation}
  \theta = \frac{T}{\pi}\,,\label{eq:15}\quad\quad\pi = \left(\frac{p}{p_g}\right)^{\frac{R}{c_p}}\,. 
\end{equation}
Variable $\theta$ is of interest in atmospheric studies because it is conserved in adiabatic processes. 
If one replaces equation~(\ref{eq:6}) with
\begin{equation}
\frac{\partial (\rho \theta)}{\partial t} + \nabla \cdot (\rho  \mathbf{u} \theta)  - \frac{\mu_a}{Pr} \Delta \theta = 0, \label{eq:16} 
\end{equation}
variable $\theta$ is computed directly, rather than obtained from post-processing of the computed variables. As a consequence of thus replacement,  the thermodynamics equations of state \eqref{eq:4} is recast into:
\begin{equation}
    p = p_g\left(\frac{\rho R \theta}{p_g}\right)^{\frac{c_p}{c_v}}. \label{eq:17}
\end{equation}
Similar to the splitting adopted for the pressure in \eqref{eq:p_splitCE2} and the density in \eqref{eq:rhoPrime} , we split the potential temperature 
into a value in hydrostatic balance $\theta_0$ and fluctuation $\theta^ \prime$ over it, i.e.,
\begin{align} 
\theta =  \theta_0 + \theta ^ \prime.
\label{pot_temp2}
\end{align}
In summary, the Euler equations considered in this work are given by
(\ref{eq:1}),\eqref{eq:9}-(\ref{eq:13}),(\ref{eq:16}),(\ref{eq:17}). 

\subsection{Discretization}\label{sec:splitting}
We adopt the backward Euler scheme for the discretization of the time derivatives in problem (\ref{eq:1}), \eqref{eq:9}-(\ref{eq:13}), (\ref{eq:16}), (\ref{eq:17}) and obtain the following time-discrete formulation:
given $\rho^0$, $\boldsymbol{u}^0$, $\theta^0$ and $p^0$,  for $n\geq0$ find $\rho^{n+1}$, $\rho^{\prime,n+1}$, $\boldsymbol{u}^{n+1}$, $\theta^{n+1}$, $p^{n+1}$, and $p^{\prime,n+1}$ such that:
\begin{align}
    & \frac{\rho^{n+1}}{\Delta t}\,+\nabla \cdot\left(\rho^{n+1}\boldsymbol{u}^{n+1}\right) = b^{n+1}_{\rho} \label{eq:18}\\[2mm]
    & \frac{\rho^{n+1} \mathbf{u}^{n+1}}{\Delta t} + \nabla \cdot \big( \rho^{n+1} \mathbf{u}^{n+1} \otimes \mathbf{u}^{n+1} \big) + \nabla p^{\prime,n+1} + \rho^{\prime , n+1} g \boldsymbol{\hat{k}} - \mu_a \Delta \mathbf{u}^{n+1} = \mathbf{b}^{n+1}_{\boldsymbol{u}}\label{eq:21}\\[2mm]
    & \frac{\rho^{n+1} \theta^{n+1}}{\Delta t} + \nabla \cdot \big( \rho^{n+1} \mathbf{u}^{n+1} \theta^{n+1} \big) - \frac{\mu_a}{Pr} \Delta \theta^{n+1} = b_{\theta}^{n+1}\label{eq:22}\\[2mm]
        & p^{n+1} = p_g \left( \frac{\rho^{n+1} R \theta^{n+1}}{p_g} \right)^{ \frac{c_p}{c_v}}, \label{eq:23}\\[2mm]
    & p^{\prime, n+1} = p^{n+1} \,-p_0 \label{eq:19} \\[2mm]
    & \rho^{\prime,n+1} = \rho^{n+1} - \rho_0 \label{eq:20}
\end{align}
where $b_{\rho}^{n+1} = \rho^n / \Delta t$, $\boldsymbol{b}_{\boldsymbol{u}}^{n+1} = \rho^n \boldsymbol{u}^n/\Delta t $ and $b_{\theta}^{n+1} = \rho^n\theta^n / \Delta t$. 
{We emphasize the distinction between the subscript and superscript notation: while $p_0$ and $\rho_0$ denote the hydrostatic pressure and density obtained from \eqref{eq:9}-\eqref{eq:14},  $p^0$ and $\rho^0$ represent the pressure and density at the initial time $t^0$. We further note that the pressure and density fluctuations are initialized according to \eqref{eq:p_splitCE2} and \eqref{eq:rhoPrime}, respectively, using the initial values, $p^0$ and $\rho^0$,  together with the equilibrium fields $p_0$ and $\rho_0$.} For more details, see, e.g., \cite{GIRFOGLIO2025106510}.

\vspace*{0.2cm}

To contain the computational cost, 
we resort to a splitting approach which consists of three steps. \\
Given $\rho^0$, $\mathbf{u}^0$, $\theta^0$, and $p^0$ for $n \geq 0$, perform the following steps:
\begin{itemize}
    \item[-] \textbf{Step 1:} find the first intermediate density $\rho^{n+\frac{1}{3}}$, the intermediate density fluctuation $\rho'^{,n+\frac{1}{3}}$ and the intermediate velocity $\mathbf{u}^{n+\frac{1}{3}}$ such that: 
\begin{align}  
&\frac{\rho^{n+\frac{1}{3}} }{\Delta t} + \nabla \cdot \left( \rho^n \mathbf{u}^n \right) = b_{\rho}^{n+1},\label{eq:24}\\[2mm]
&\frac{\rho^{n+\frac{1}{3}} \mathbf{u}^{n+\frac{1}{3}} }{\Delta t} + \nabla \cdot \left( \rho^{n} \mathbf{u}^n \otimes \mathbf{u}^{n+\frac{1}{3}} \right) 
+ \nabla p'^{,n} + \rho'^{,n+\frac{1}{3}} g \boldsymbol{\hat{k}}  - \mu_a\Delta \mathbf{u}^{n+\frac{1}{3}} = \mathbf{b}_{\mathbf{u}}^{n+1},\label{eq:26} \\[2mm]
&\rho'^{,n+\frac{1}{3}} = \rho^{n+\frac{1}{3}} - \rho_0.\label{eq:25}
\end{align}
    \item[-] \textbf{Step 2:} find the potential temperature $\theta^{n+1}$ and the second intermediate density $\rho^{n+\frac{2}{3}}$ such that:
\begin{align}
&\frac{\rho^{n+\frac{1}{3}} \theta^{n+1} }{\Delta t} + \nabla \cdot \left( \rho^{n} \mathbf{u}^n \theta^{n+1} \right) - \frac{\mu_a}{Pr} \Delta \theta^{n+1} = b_{\theta}^{n+1},\quad\quad\quad\quad\quad\quad\quad\quad \textcolor{white}{frase}\label{eq:27}\\
&p^n = p_g \left( \frac{\rho^{n+\frac{2}{3}} R \theta^{n+1}}{p_g} \right)^{ \frac{c_p}{c_v}}. \label{eq:28} 
\end{align}
\item[-] \textbf{Step 3:} find the end-of-step velocity $\mathbf{u}^{n+1}$, density $\rho^{n+1}$, pressure $p^{n+1}$ and pressure fluctuation $p'^{,n+1}$ such that:  
\begin{align}  
&\frac{\rho^{n+\frac{1}{3}} \mathbf{u}^{n+1}}{\Delta t} + \nabla \cdot \left( \rho^{n} \mathbf{u}^n \otimes \mathbf{u}^{n+1} \right) 
+ \nabla p'^{,n+1} + \rho'^{,n+\frac{1}{3}} g \boldsymbol{\hat{k}} -\mu_a \Delta \mathbf{u}^{n+\frac{1}{3}} = \mathbf{b}_{\mathbf{u}}^{n+1},\label{eq:29} \\[2mm]
&\frac{\rho^{n+1} }{\Delta t} + \nabla \cdot \left( \rho^{n+\frac{2}{3}} \mathbf{u}^{n+1} \right) = b_{\rho}^{n+1},\label{eq:32} \\[2mm]
&p^{n+1} = p_g \left( \frac{\rho^{n+1} R \theta^{n+1}}{p_g} \right)^{ \frac{c_p}{c_v}},\label{eq:31}\\[2mm]
&p'^{,n+1} =p^{n+1} - p_0 .\label{eq:30}
\end{align}
\end{itemize}
The equations in Steps 1-3 are 
approximated in space using 
a finite volume method. 
We refer to \cite{GQR_OF_clima,GIRFOGLIO2025106510} for additional details about this discretization
and the treatment of the hydrostatic balance. 
In \cite{GIRFOGLIO2025106510}, one can also find a thorough assessment of the formulation and the numerical methodology through benchmark tests involving dry air flow over a flat terrain and orography.

\section{Numerical results}\label{sec:flatBench}
In this section, the proposed IREE-driven mesh adaptation strategy is evaluated and compared with the PMA algorithm using two well-established two-dimensional benchmarks for atmospheric flows: the rising thermal bubble \cite{Marras2015,Ahmad2007,Ahmad2018,Feng2021} and the density current \cite{Marras2015,Ahmad2007,Carpenter1990,Giraldo2008,Marras2013,Straka1993}. 
Both test cases involve a perturbation of a neutrally stratified atmosphere with uniform background potential temperature over a flat terrain. 
Neither of these test cases has an exact solution, hence  comparisons can only be made relative to other numerical data available in the literature or numerical solutions associated with a very fine computational grid. Out of the several variations of these benchmarks, we select the settings from 
\cite{Ahmad2007,Feng2021} for the rising thermal bubble and from \cite{Carpenter1990,Straka1993}
for the density current. 
For the density current, we also present a novel three-dimensional version of the benchmark.

The IREE and PMA algorithms for atmospheric flows are implemented in GEA (Geophysical and Environmental Applications) \cite{GEA,GirfoglioFVCA10, GQR_OF_clima,CGQR}, an open-source package for atmosphere and ocean modeling based on the finite volume C++ library OpenFOAM\textsuperscript{\textregistered}. OpenFOAM\textsuperscript{\textregistered} manages octree structures and extends its native mesh adaptation algorithm (i.e., the PMA algorithm in Sec.~\ref{sec:plain}) to quadtree structures through
a well established external library \cite{Darbar2022}. In particular, we drive mesh adaptation in both the IREE and PMA algorithms by means of the gradient of the potential temperature, $\mathbf{Q}=\nabla\theta$.

\subsection{Two-dimensional rising thermal bubble}\label{sec:2Dbubble}

The computational domain in the $xz$-plane is $\Omega = [0,5000] \times [0,10000] \, \mathrm{m^2}$ and the time interval of interest is $(0, 1020]$ s. Impenetrable, free-slip boundary conditions are imposed on all domain walls. The initial potential temperature is
\begin{equation}
\theta^0 = \left\{
\begin{array}{ll}
300 + 2 \left[ 1 - \frac{r}{r_0} \right] \qquad &\text{if } r \leq r_0 = 2000 \, \mathrm{m}\\[2mm] 300 &\text{ otherwise},
\end{array}
\right.
\end{equation}
with radius $r = \sqrt{(x - x_c)^2 + (z - z_c)^2}$ and centre $(x_c, z_c) = (0, 2000) \, \mathrm{m}$. 
This initial condition refers to a neutrally stratified atmosphere, with uniform background potential temperature of $300 \, \mathrm{K}$ perturbed by a circular bubble of warmer air. The initial density $\rho^0$ is computed by solving the hydrostatic balance \eqref{eq:14}, while the initial velocity field is zero everywhere. In addition, following \cite{Ahmad2007}, we set $\mu_a = 15$ and $Pr = 1$ in \eqref{eq:21} and \eqref{eq:22}. 

We consider three different initial meshes with uniform resolution $h = \Delta x = \Delta z =$ 125, 62.5, 15.625 m. 
The time step is set to $\Delta t = 0.1$ s for all the simulations. The two coarser meshes are used to assess the performance of the IREE algorithm against both the PMA algorithm and a reference solution computed on the finest mesh without adaptation.
We note that 
it is rather common to consider the results obtained with mesh $h = 15.625$ m as the reference when $\mu_a = 15$ \cite{Ahmad2007,Marras2015,GIRFOGLIO2025106510,GQR_OF_clima,CGQR,Clinco2024}.

With the purpose of comparing the IREE and PMA algorithms, for both 
we set $N_\text{max}=2\times10^6$ and 
$I_r = 5$. 
All the other parameter values for the two algorithms are reported in Tab.~\ref{tab:1}.
These values were obtained through a trial-and-error procedure for this specific benchmark and may require retuning for different test cases.
%
\begin{table}[htb!]
    \centering
    \begin{minipage}{0.45\textwidth}  
        \centering
        \begin{tabular}{|c|c|c|c|c|}
            \hline
            \multicolumn{5}{|c|}{\textbf{IREE}} \\ \hline
            \textbf{Mesh size} & $\delta_1$ & $\delta_2$ & $tol$ & $l_\text{max}$\\ \hline
            62.5 & 6.0 & 0.90 & 1e-3 & 1  \\ \hline
            125 & 1.5 & 0.50 & 1e-3 & 2  \\ \hline
        \end{tabular}
    \end{minipage}
    \begin{minipage}{0.45\textwidth}  
        \centering
        \begin{tabular}{|c|c|c|c|c|}
            \hline
            \multicolumn{5}{|c|}{\textbf{PMA}} \\ \hline
            \textbf{Mesh size} & $\alpha_\text{min}^\text{ref}$ & $\alpha_\text{max}^\text{ref}$ & $l_\text{max}$ & $l_c$ \\ \hline
            62.5 & 0.1 & 0.9 & 1 & 10  \\ \hline
            125 & 0.1 & 0.9 & 2 & 10  \\ \hline
        \end{tabular}
    \end{minipage}%
        \caption{2D rising thermal bubble: parameter values for the IREE and PMA algorithms.}\label{tab:1}
\end{table}

\noindent
Figure~\ref{fig:BubbleInTime2D_125} shows the evolution of the computed potential temperature and the corresponding
mesh adapted by the IREE and PMA algorithms,
starting from the mesh size $h=125$ m. 
IREE algorithm refines
the mesh over a larger region than the PMA approach. {Next, we  investigate whether this additional refinement translates into a more accurate solution.}
\begin{figure}[htb!]
    \centering
    
    \begin{overpic}[width=0.24\textwidth]{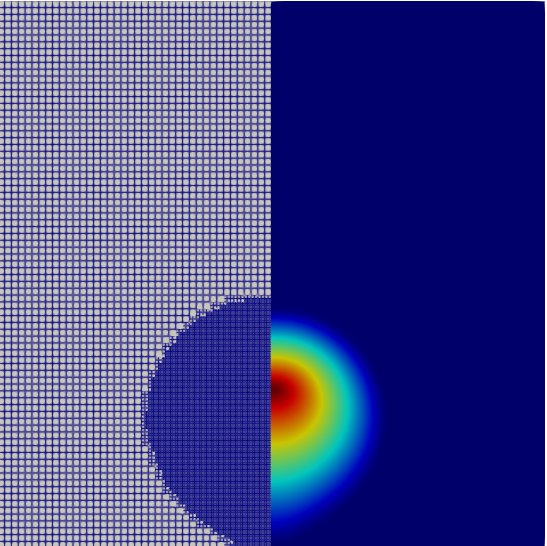}
        \put(53,80){\textcolor{white}{\footnotesize{\shortstack{$t=250$ s}}}}
        \put(53,88){\textcolor{white}{\footnotesize{\shortstack{IREE}}}}
    \end{overpic}
    \begin{overpic}[width=0.24\textwidth]{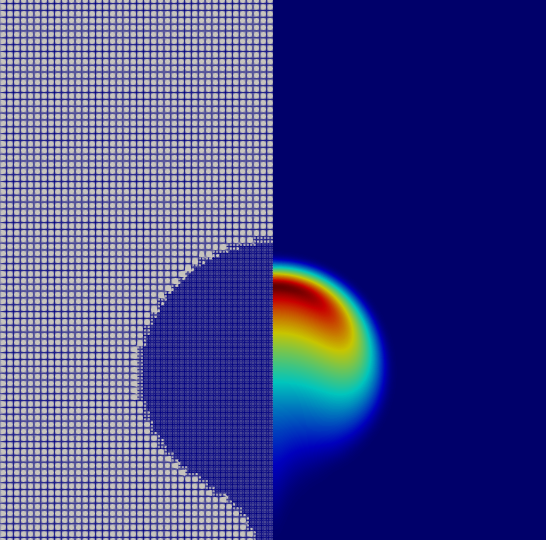}
        \put(53,80){\textcolor{white}{\footnotesize{\shortstack{$t=500$ s}}}}
        \put(53,88){\textcolor{white}{\footnotesize{\shortstack{IREE}}}}
    \end{overpic}
    \begin{overpic}[width=0.24\textwidth]{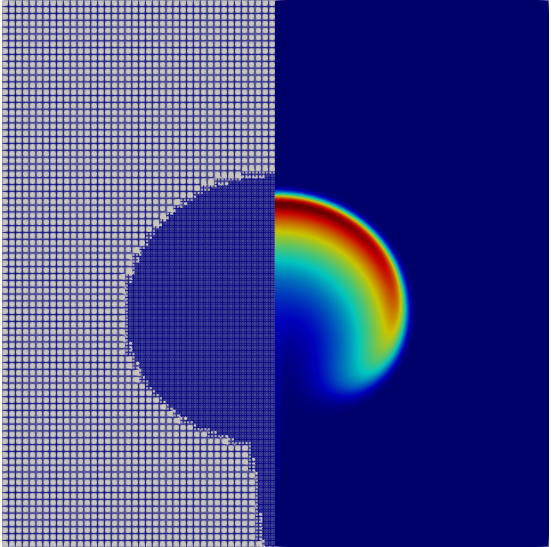}
        \put(53,80){\textcolor{white}{\footnotesize{\shortstack{$t=750$ s}}}}
        \put(53,88){\textcolor{white}{\footnotesize{\shortstack{IREE}}}}
    \end{overpic}
    \begin{overpic}[width=0.24\textwidth]{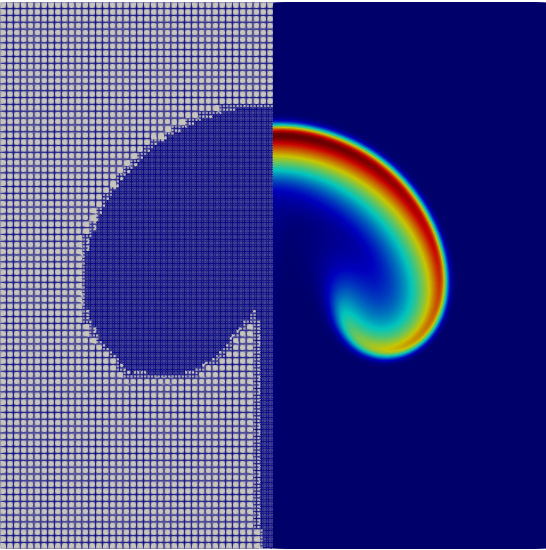}
        \put(53,80){\textcolor{white}{\footnotesize{\shortstack{$t=1020$ s}}}}
        \put(53,88){\textcolor{white}{\footnotesize{\shortstack{IREE}}}}
    \end{overpic}
\\
\vskip .1cm
    \begin{overpic}[width=0.24\textwidth]{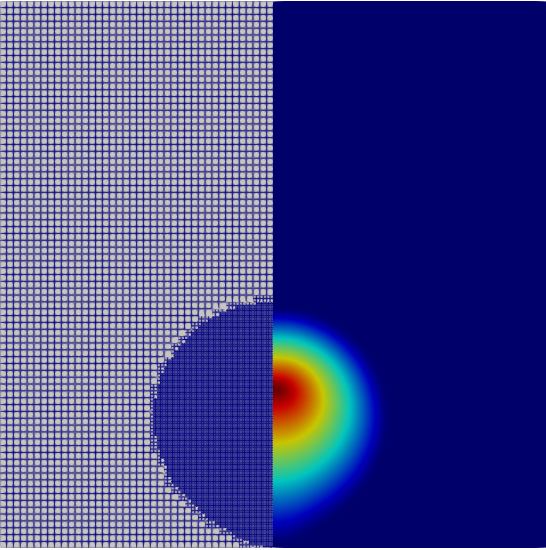}
        \put(53,80){\textcolor{white}{\footnotesize{\shortstack{$t=250$ s}}}}
        \put(53,88){\textcolor{white}{\footnotesize{\shortstack{PMA}}}}
    \end{overpic}
        \begin{overpic}[width=0.24\textwidth]{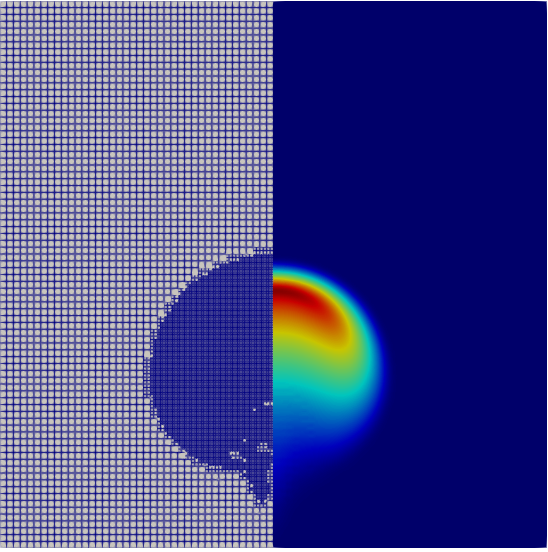}
        \put(53,80){\textcolor{white}{\footnotesize{\shortstack{$t=500$ s}}}}
        \put(53,88){\textcolor{white}{\footnotesize{\shortstack{PMA}}}}
    \end{overpic}
    \begin{overpic}[width=0.24\textwidth]{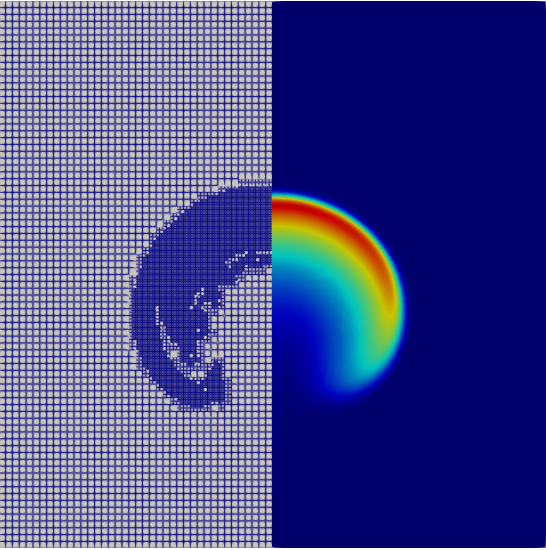}
        \put(53,80){\textcolor{white}{\footnotesize{\shortstack{$t=750$ s}}}}
        \put(53,88){\textcolor{white}{\footnotesize{\shortstack{PMA}}}}
    \end{overpic}
    \begin{overpic}[width=0.24\textwidth]{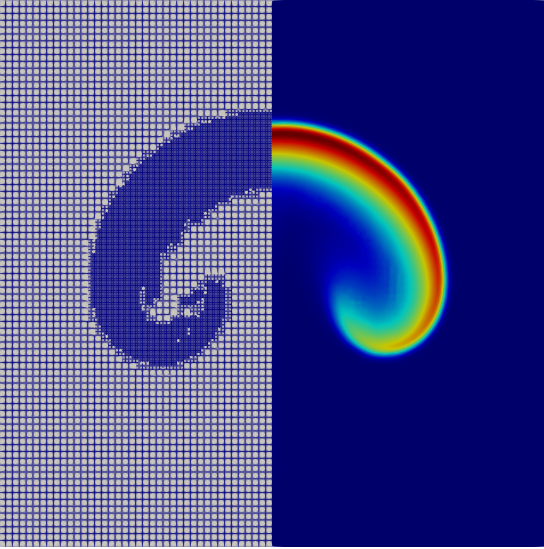}
        \put(53,80){\textcolor{white}{\footnotesize{\shortstack{$t=1020$ s}}}}
        \put(53,88){\textcolor{white}{\footnotesize{\shortstack{PMA}}}}
    \end{overpic}
    \vspace{5pt}
    \begin{overpic}[width=0.50\textwidth]{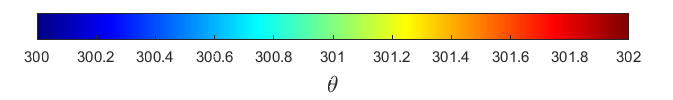}
    \end{overpic}
    \caption{2D rising thermal bubble: 
    evolution of the potential temperature and associated mesh adapted by IREE (top) and PMA (bottom) algorithms starting from a mesh size $h=125$ m.}
    \label{fig:BubbleInTime2D_125}
\end{figure}
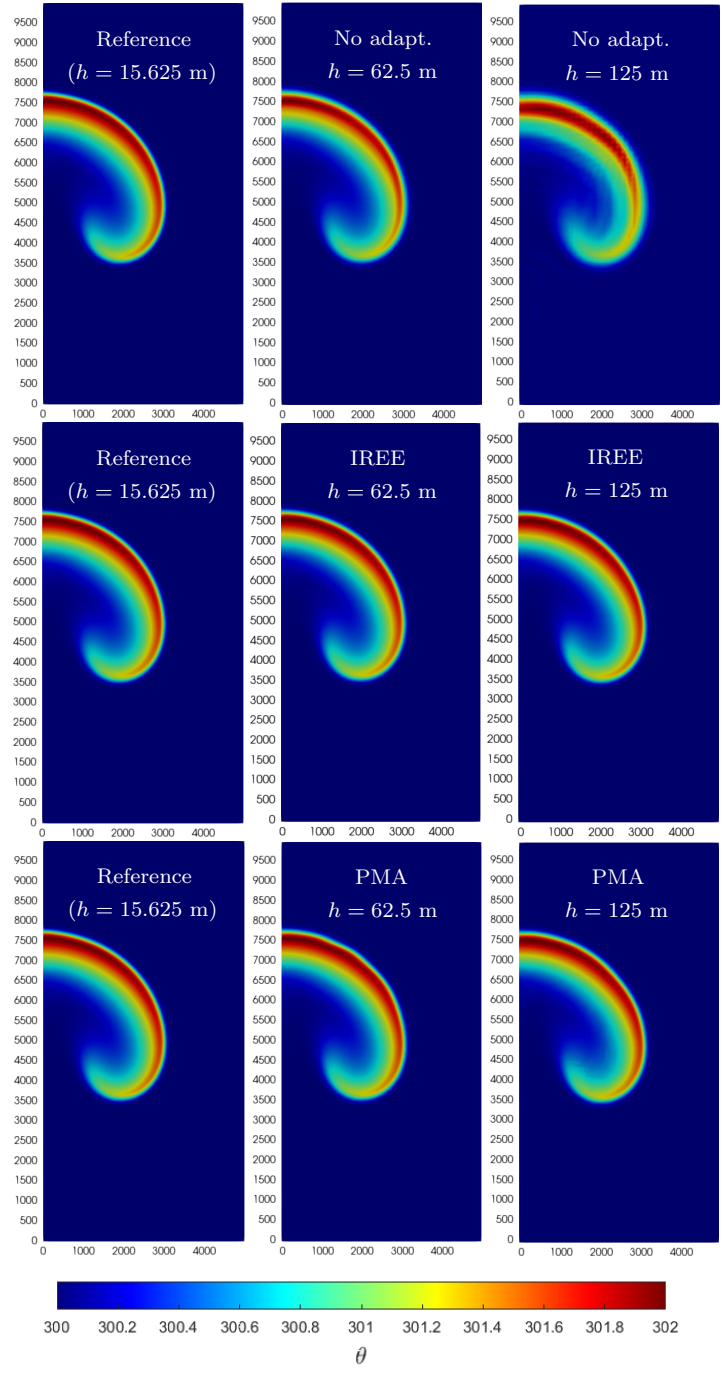
\begin{figure}[htb!]
    \centering
    
    \begin{overpic}[width=0.20\textwidth, grid=false]{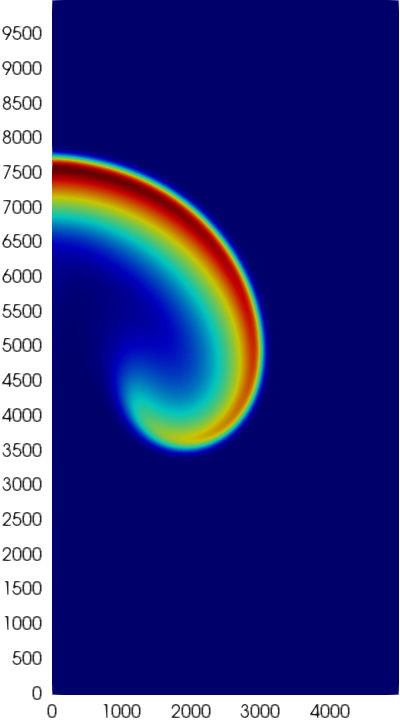}
        \put(20,90){\textcolor{white}{\footnotesize{\shortstack{ Reference}}}}
        \put(13,82){\textcolor{white}{\footnotesize{\shortstack{ ($h = 15.625$ m)}}}}
    \end{overpic}
    \begin{overpic}[width=0.20\textwidth, grid = false]{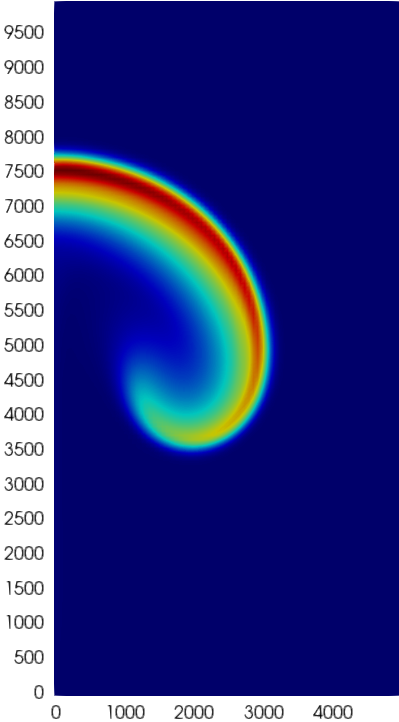}
        \put(20,90){\textcolor{white}{\footnotesize{\shortstack{No adapt.}}}}
        \put(18.5,82){\textcolor{white}{\footnotesize{\shortstack{ $h = 62.5$ m}}}}
    \end{overpic}
    \begin{overpic}[width=0.20\textwidth]{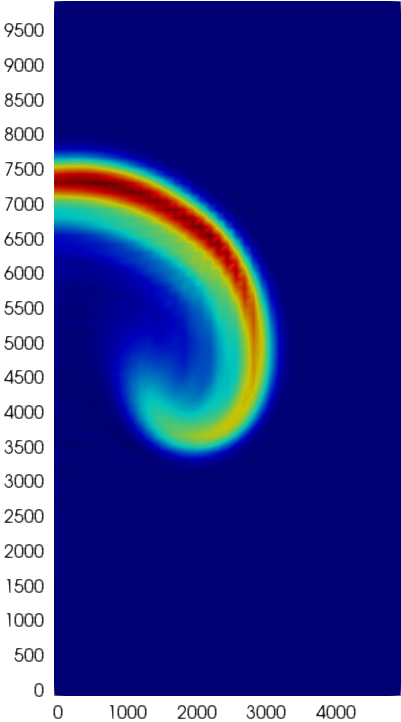}
        \put(20,90){\textcolor{white}{\footnotesize{\shortstack{No adapt.}}}}
        \put(18.5,82){\textcolor{white}{\footnotesize{\shortstack{ $h = 125$ m}}}}
    \end{overpic}
    \\
    \begin{overpic}[width=0.20\textwidth]{Images/RB_contour/RB_15_625m.png}
        \put(20,90){\textcolor{white}{\footnotesize{\shortstack{ Reference}}}}
        \put(13,82){\textcolor{white}{\footnotesize{\shortstack{ ($h = 15.625$ m)}}}}
        \end{overpic}
    \begin{overpic}[width=0.20\textwidth]{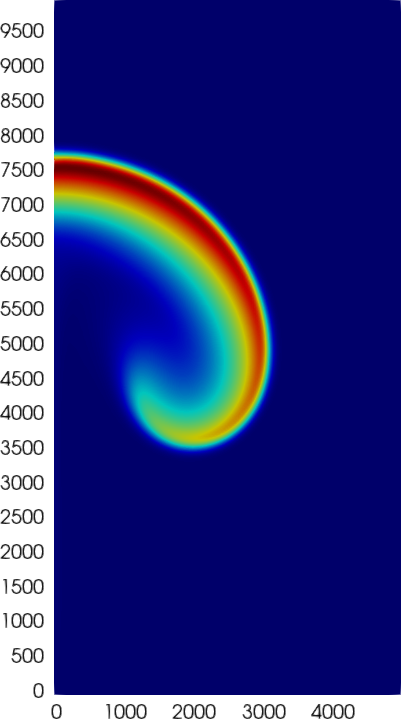}
        \put(24,90){\textcolor{white}{\footnotesize{\shortstack{IREE}}}}
        \put(18.5,82){\textcolor{white}{\footnotesize{\shortstack{ $h = 62.5$ m}}}}
    \end{overpic}
    \begin{overpic}[width=0.20\textwidth]{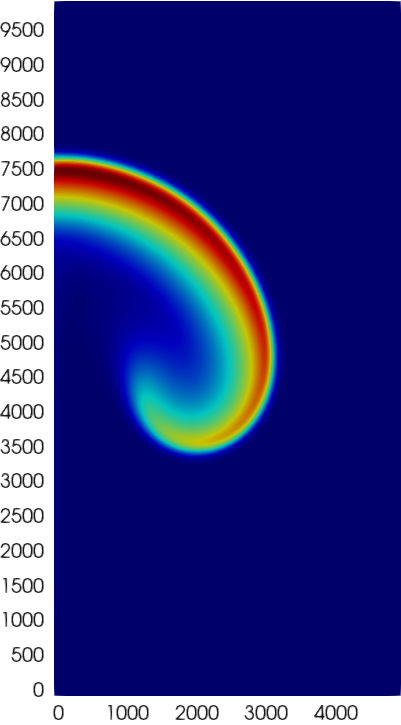}
        \put(24,90){\textcolor{white}{\footnotesize{\shortstack{IREE}}}}
        \put(18.5,82){\textcolor{white}{\footnotesize{\shortstack{ $h = 125$ m}}}}
    \end{overpic}
    \\
    \begin{overpic}[width=0.20\textwidth]{Images/RB_contour/RB_15_625m.png}
                \put(20,90){\textcolor{white}{\footnotesize{\shortstack{ Reference}}}}
        \put(13,82){\textcolor{white}{\footnotesize{\shortstack{ ($h = 15.625$ m)}}}}
    \end{overpic}
    \begin{overpic}[width=0.20\textwidth]{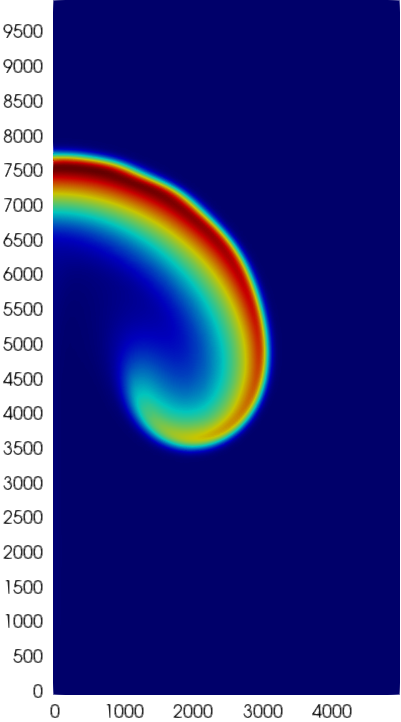}
        \put(25,90){\textcolor{white}{\footnotesize{\shortstack{PMA}}}}
        \put(18.5,82){\textcolor{white}{\footnotesize{\shortstack{ $h = 62.5$ m}}}}
    \end{overpic}
    \begin{overpic}[width=0.20\textwidth]{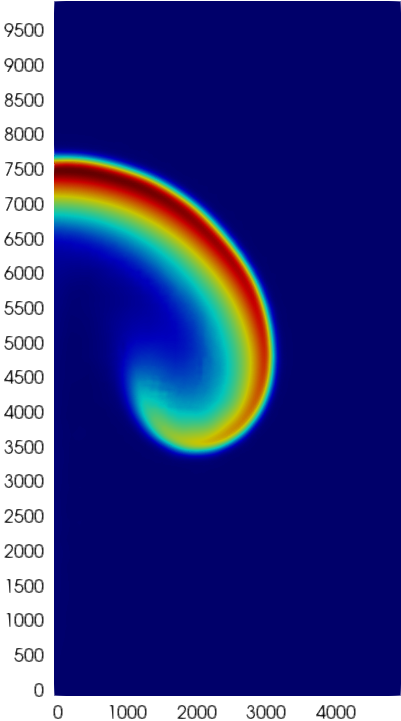}
        \put(25,90){\textcolor{white}{\footnotesize{\shortstack{PMA}}}}
        \put(18.5,82){\textcolor{white}{\footnotesize{\shortstack{ $h = 125$ m}}}}
    \end{overpic}

\vskip .2cm
        \begin{overpic}[width=0.60\textwidth]{Images/colorbar/colorbar.png}
        \put(53,85){\textcolor{white}{\footnotesize{}}}
    \end{overpic}
    
    \caption{2D rising thermal bubble:  potential temperature at time $t=1020$ s provided by the reference simulation (first column) compared against the solution computed with mesh size $h = 62.5$ m (second column) and $h = 125$ m (third column), with no mesh adaptation (first row) and with mesh adaptation driven by the IREE (second row) and PMA (third row) algorithm.}
    \label{fig:BubbleContour2D}
\end{figure}

\begin{figure}[htb!]
    \centering
    \begin{overpic}[width=0.20\textwidth]{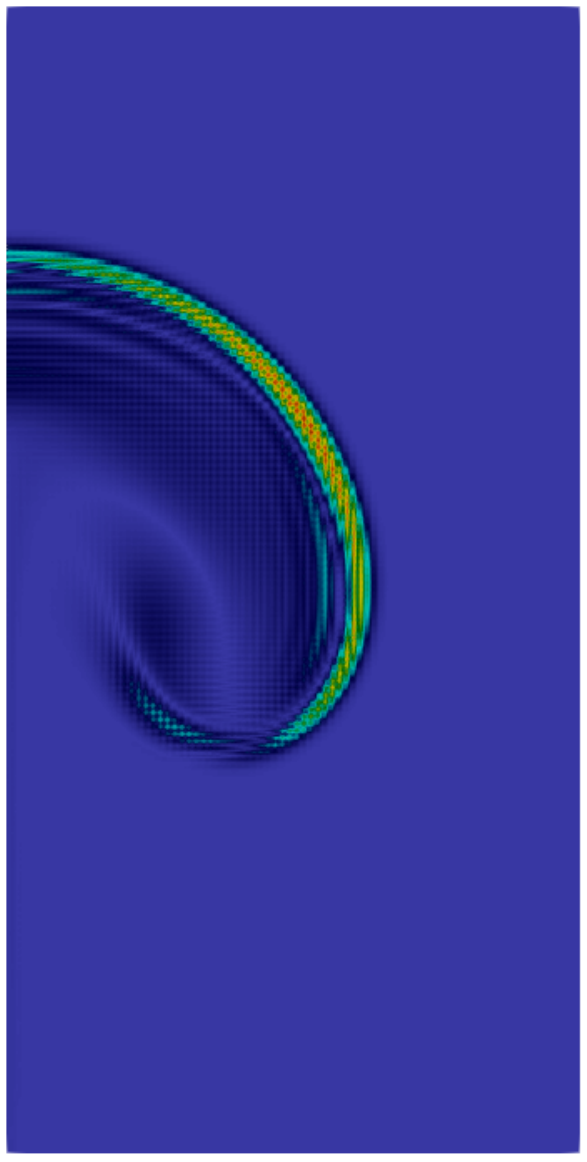}
        \put(15,90){\textcolor{white}{\footnotesize{\shortstack{ No adapt.}}}}
        \put(13,82){\textcolor{white}{\footnotesize{\shortstack{ $h = 62.5$ m}}}}
        \end{overpic}
    \begin{overpic}[width=0.20\textwidth]{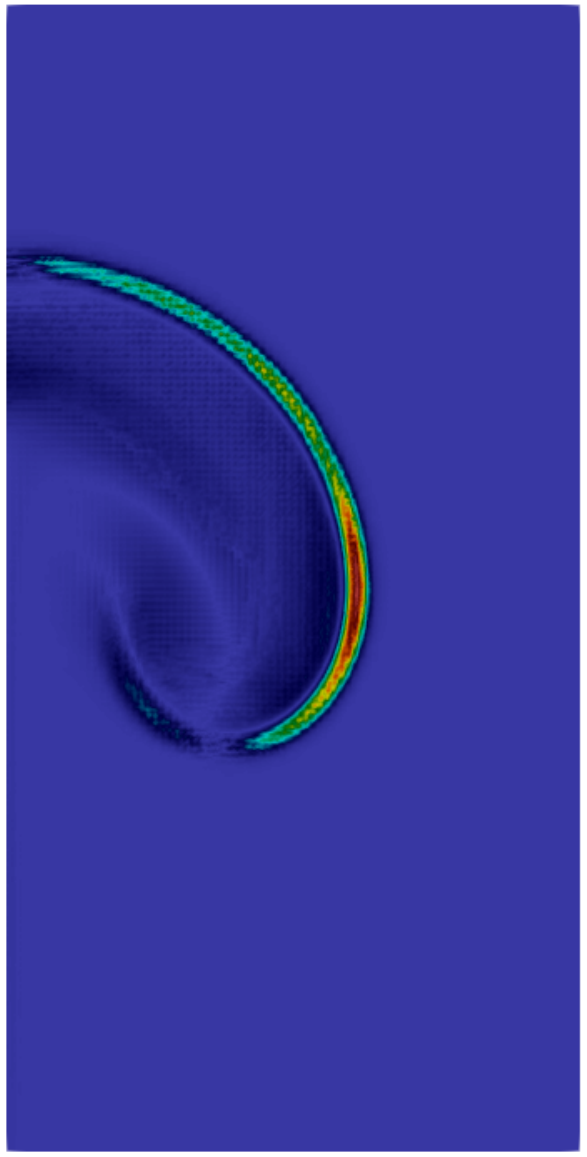}
        \put(19,90){\textcolor{white}{\footnotesize{\shortstack{PMA}}}}
        \put(14,82){\textcolor{white}{\footnotesize{\shortstack{ $h = 62.5$ m}}}}
    \end{overpic}
    \begin{overpic}[width=0.20\textwidth]{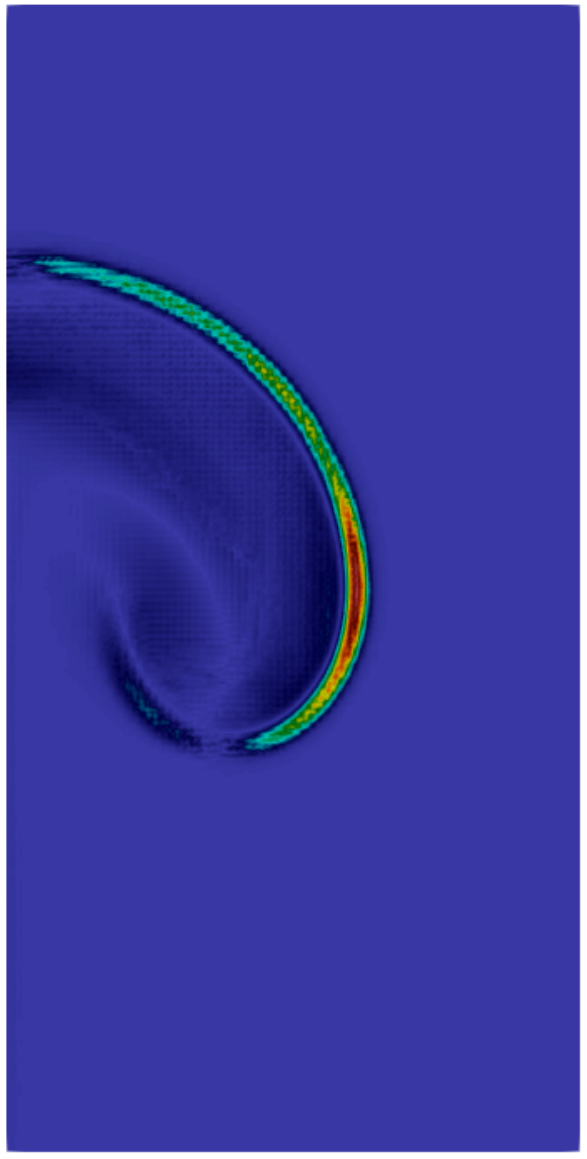}
        \put(18,90){\textcolor{white}{\footnotesize{\shortstack{IREE}}}}
        \put(13,82){\textcolor{white}{\footnotesize{\shortstack{ $h = 62.5$ m}}}}
        \end{overpic}
    \\
        \begin{overpic}[width=0.20\textwidth]{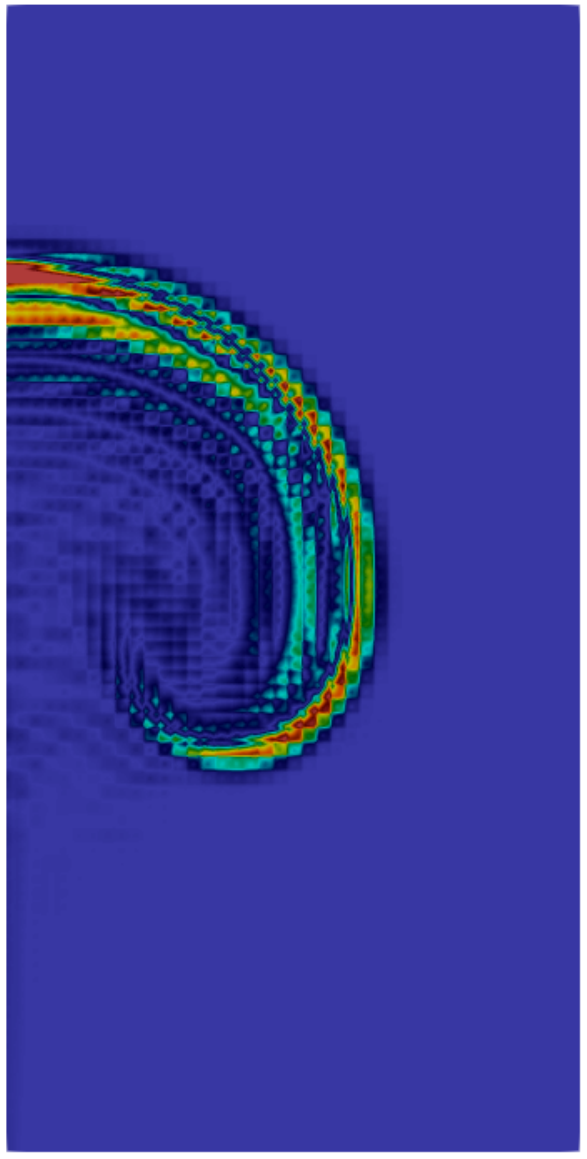}
        \put(15,90){\textcolor{white}{\footnotesize{\shortstack{No adapt.}}}}
        \put(13,82){\textcolor{white}{\footnotesize{\shortstack{ $h = 125$ m}}}}
        \end{overpic}
    \begin{overpic}[width=0.20\textwidth]{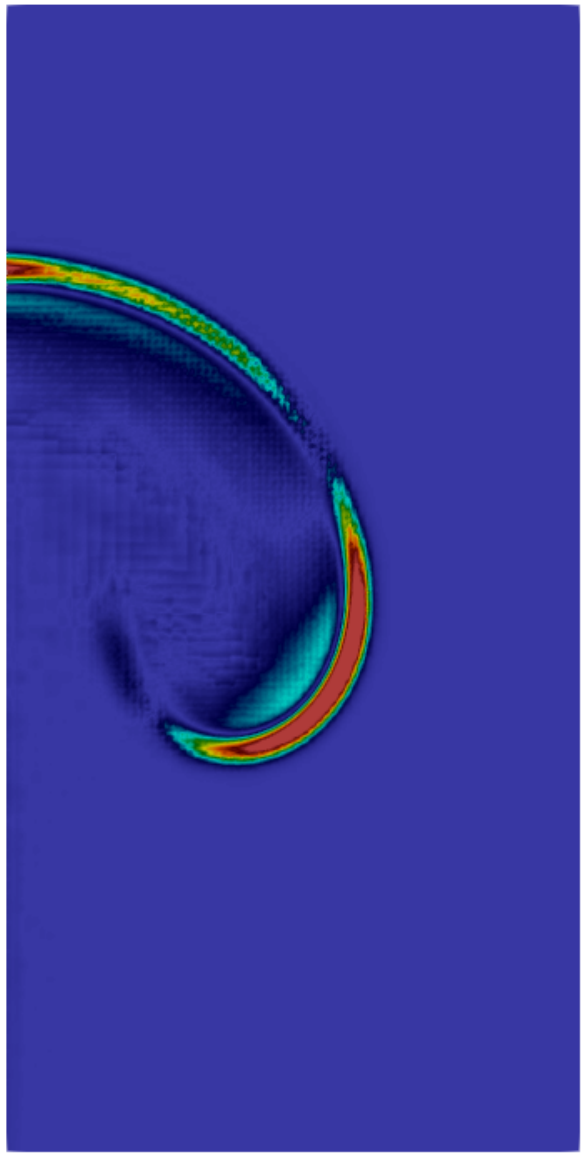}
        \put(19,90){\textcolor{white}{\footnotesize{\shortstack{PMA}}}}
        \put(14,82){\textcolor{white}{\footnotesize{\shortstack{ $h = 125$ m}}}}
    \end{overpic}
    \begin{overpic}[width=0.20\textwidth]{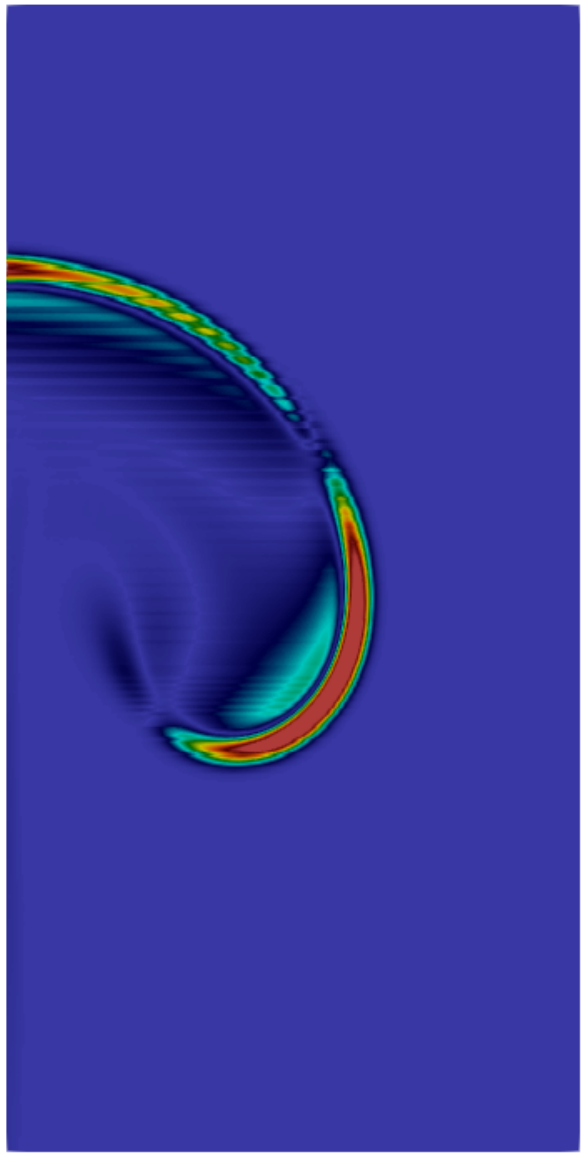}
        \put(18,90){\textcolor{white}{\footnotesize{\shortstack{IREE}}}}
        \put(14,82){\textcolor{white}{\footnotesize{\shortstack{ $h = 125$ m}}}}
    \end{overpic}

\vskip .2cm
        \begin{overpic}[width=0.60\textwidth]{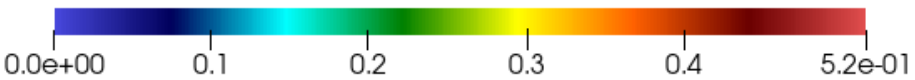}
        \put(53,85){\textcolor{white}{\footnotesize{}}}
    \end{overpic}
    
    \caption{2D rising thermal bubble: absolute error for the potential temperature at time $t=1020$ s for the simulations with no mesh adaptation (first column), PMA (second column) and IREE (third column) algorithms with respect to the reference simulation (mesh $h = 15.625$ m) starting from meshes $h = 62.5$ m (first row) and $h = 125$ m (second row).}
    \label{fig:BubbleContour2D_error}
\end{figure}
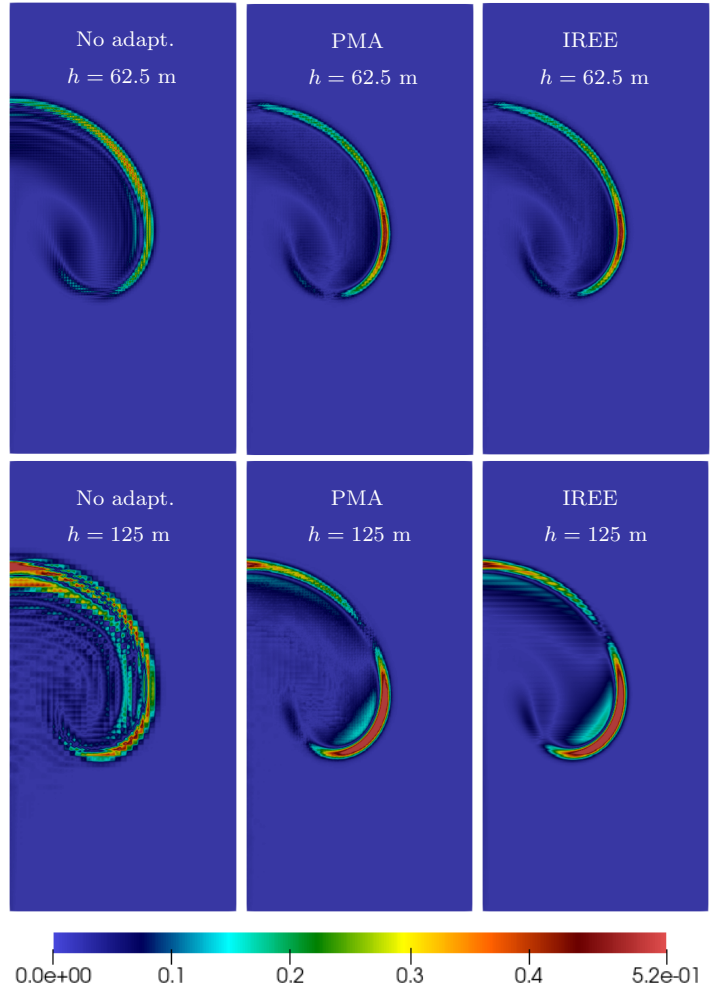

\begin{figure}[htb!]
    \centering  
    \begin{subfigure}{0.48\textwidth}
        \centering
        \includegraphics[width=\textwidth]{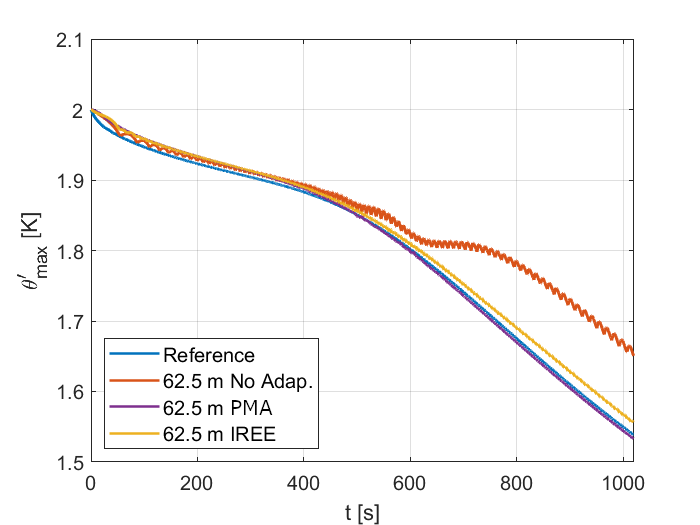}
        \label{fig:immagine1}
    \end{subfigure}
    \begin{subfigure}{0.5\textwidth}
        \centering
        \includegraphics[width=\textwidth]{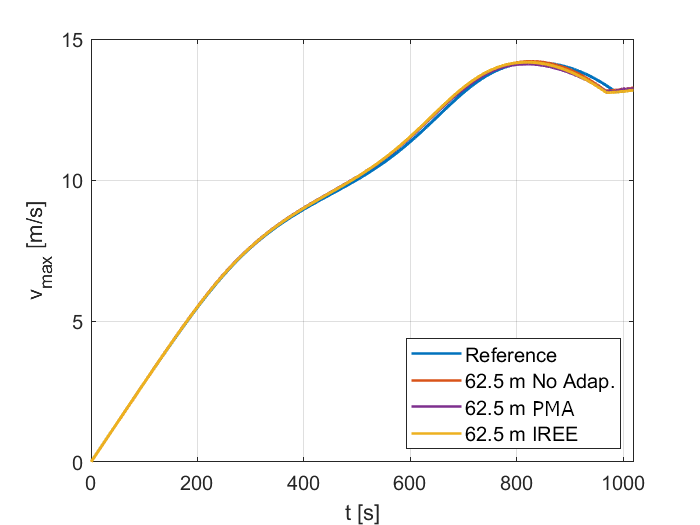}

        \label{fig:immagine2}
    \end{subfigure}
    
    \begin{subfigure}{0.48\textwidth}
        \centering
        \includegraphics[width=\textwidth]{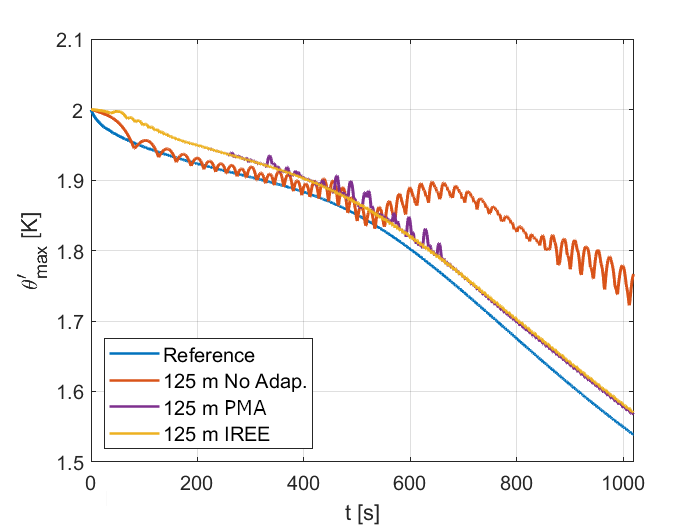}
        \label{fig:immagine1b}
    \end{subfigure}
    \begin{subfigure}{0.48\textwidth}
        \centering
        \includegraphics[width=\textwidth]{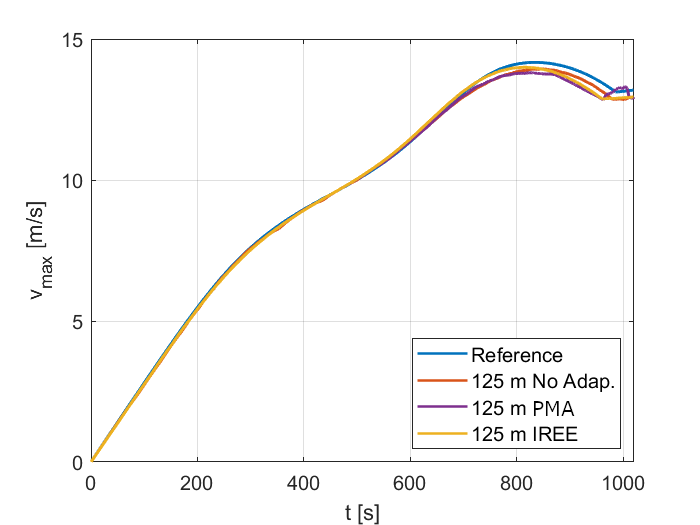}
        \label{fig:immagine2b}
    \end{subfigure}

    \caption{
    2D rising thermal bubble: evolution of the maximum potential temperature fluctuation, $\theta'_{max}$, (left) and of the maximum vertical component of the velocity, $w_{max}$, (right) computed with a mesh size $h = 62.5$ m (top) and $h = 125$ m (bottom). The reference values are associated with a fixed uniformly fine grid with mesh size $h = 15.625$ m.}
    \label{fig:BubbleMatlabPlot2D}
\end{figure}

Fig.~\ref{fig:BubbleContour2D} compares the potential temperature
at $t=1020$ s (i.e., at the end of the time interval under consideration) computed without and with mesh adaptation by the IREE or PMA algorithms. From the top row, 
which shows the results computed without mesh adaptation, 
we see a degradation in accuracy as the mesh gets coarser. The results associated with the mesh sizes
$h = 65.2, 125$ m become more accurate when using a mesh adaptation algorithm, as evidenced by their improved agreement with the reference solution. However, Fig.~\ref{fig:BubbleContour2D} does not clearly reveal whether the IREE algorithm outperforms the PMA adaptation procedure. 
To obtain a more quantitative comparison, Fig.~\ref{fig:BubbleContour2D_error} displays the absolute error in the potential temperature with respect to the reference solution for all the simulations shown in Fig.~\ref{fig:BubbleContour2D}.
The checkerboard patterns visible in the error plots for the simulation without mesh adaptation indicate the onset of numerical instabilities. 
The error plots for the results given by the 
IREE and PMA algorithms starting from mesh size $h = 62.5$ m are almost indistinguishable. As for the discretization step $h = 125$ m, the error plots are still very similar, although the one associated with the PMA algorithm reveals a more evident checkerboard-like pattern, providing a first indication of numerical instabilities that will become more apparent in Fig.~\ref{fig:BubbleMatlabPlot2D}.
%

\begin{figure}[htb!]
    \centering
    \begin{subfigure}{0.48\textwidth}
        \includegraphics[width = \textwidth]{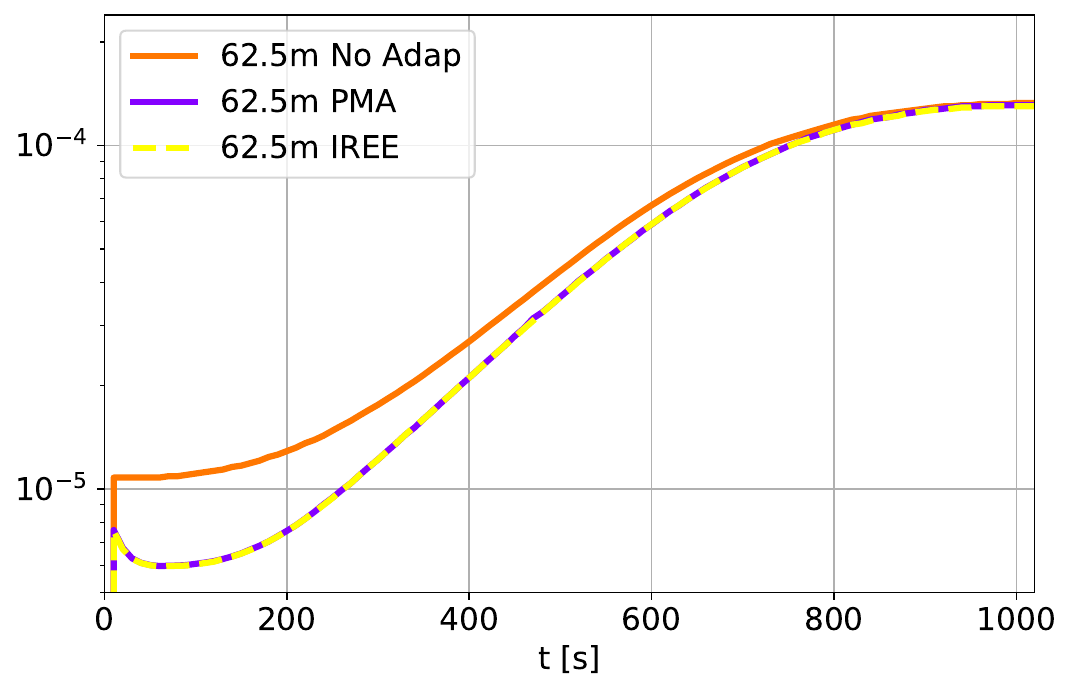}
    \end{subfigure}
    \begin{subfigure}{0.48\textwidth}
        \includegraphics[width = \textwidth]{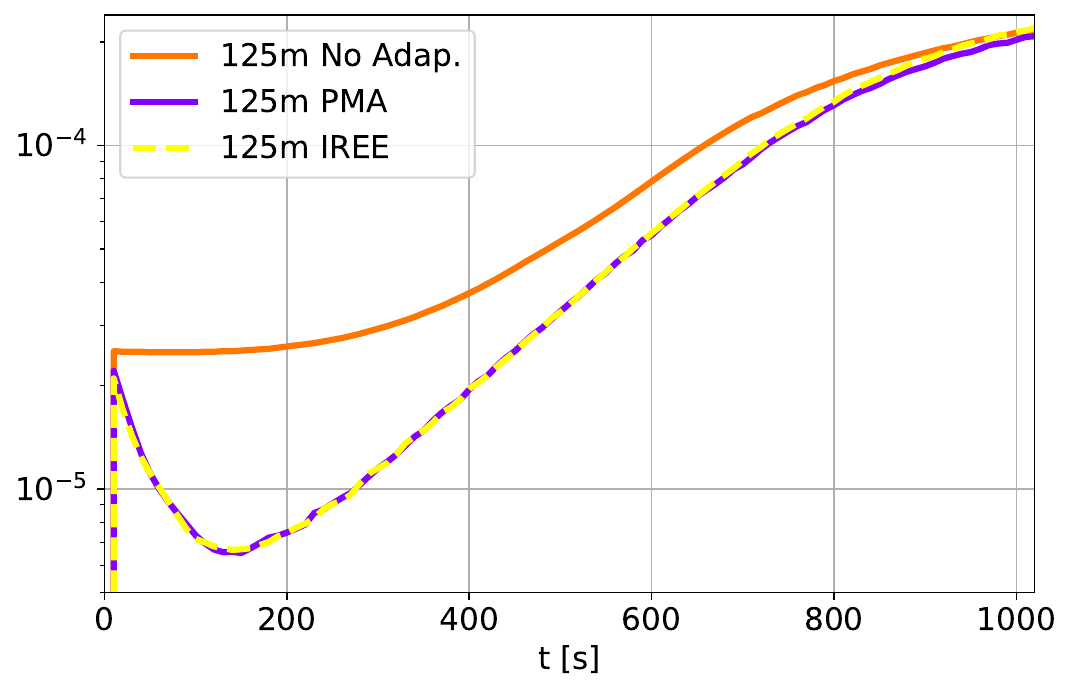}
    \end{subfigure}
    \caption{2D rising thermal bubble: the time evolution of the relative error in the $L^2(\Omega)$-norm for the potential temperature $\theta$ using mesh sizes $h = 62.5$ m (left) and $h = 125$ m (right). The reference solution in the calculation of the error is the approximation computed on a fixed uniformly fine grid with mesh size $h = 15.625$ m.}
    \label{fig:RTB_error}
\end{figure}

To further compare the performance of the IREE and PMA algorithms, we now present quantitative results.
Figure~\ref{fig:RTB_error} displays the evolution of relative error with respect to the $L^2(\Omega)$-norm, 
for all the simulations in Fig.~\ref{fig:BubbleContour2D}. 
For most of the time interval under consideration, 
the lack of mesh adaptation yields the largest error,
for both the mesh sizes, $h = 62.5, 125$ m, with the coarser mesh exhibiting the larger error, as expected.
The benefits of mesh adaptation are particularly evident during the first half of the simulation, especially for the coarser mesh with $h = 125$ m, where both the IREE and PMA algorithms achieve almost identical error reductions. 
This confirms that, for this test case, the two mesh adaptation algorithms guarantee very similar results.
In the second half of the time interval, 
the advantage of mesh adaptation in terms of relative error progressively decreases.
This should not be interpreted as a loss of accuracy due to mesh adaptation. Indeed, as already observed in Fig.~\ref{fig:BubbleContour2D_error}, the solution obtained without mesh adaptation exhibits numerical instabilities, which are largely mitigated when mesh adaptation is employed.

As a further validation step, in Fig.~\ref{fig:BubbleMatlabPlot2D} we track the evolution of two quantities of interest commonly considered in the literature \cite{Ahmad2007, Marras2013}, i.e., 
the maximum potential temperature fluctuation, $\theta^{\prime}_\text{max}$, and the maximum vertical component of the velocity, $w_\text{max}$. 

The $\theta'_{\max}$ values computed without mesh adaptation deviate significantly from the reference curve and exhibit oscillations associated with numerical instabilities. When mesh adaptation is employed, both the IREE and PMA algorithms substantially improve the agreement with the reference solution. For the finer mesh, the two adaptive strategies provide comparable results. For the coarser mesh, however, localized oscillations remain visible when using the PMA algorithm, whereas they are largely suppressed by the IREE algorithm. This suggests that the latter provides a more effective stabilization of the flow. 
As for the maximum vertical component of the velocity, 
the results indicate that this quantity is less sensitive
to variations in the mesh size and to whether or not mesh adaptation is employed. 
%

Finally, we comment on the benefits led by mesh adaptation algorithms in terms of computational performance. 
The reference simulation requires $3800$ s of computational time. Table~\ref{tab:2} reports the times associated with the IREE and PMA algorithms for mesh sizes $h = 65.2$ m and $h = 125$ m, together with the corresponding percentage time savings relative to the reference simulation.
Although both adaptive strategies yield substantial computational savings, the PMA algorithm remains the less expensive option. Figure~\ref{fig:BubbleInTime2D_125} provides a first explanation for this behavior, showing that the IREE strategy refines a larger portion of the computational domain than PMA. Moreover, the two approaches rely on different adaptation criteria: while PMA directly exploits the available quantity $\mathbf{Q}$, IREE requires the preliminary construction of a recovered approximation of the gradient field to drive mesh adaptation.
The resulting increase in computational cost may be viewed as the price to pay for a more stable numerical solution. Indeed, as shown in Fig.~\ref{fig:BubbleMatlabPlot2D}, the oscillations observed in some PMA simulations are significantly reduced when the IREE strategy is employed.
\begin{table}[htb!]
\centering
\begin{tabular}{|c||c|c|}
\hline
mesh size $h$ [m]&  IREE algorithm [s]&  PMA algorithm [s]  \\
\hline
62.5 & 761 (80\%)  & 481 (87.3\%) \\       
125 &  522 (86.3\%)   & 294 (92.3\%)\\       
\hline
\end{tabular}
\caption{2D rising thermal bubble: 
computational performance of the IREE and PMA algorithms for mesh sizes $h = 65.2, 125$ m, in terms of computational time and percentage of time saved with respect to the reference simulation ($3800$ s).}\label{tab:2}
\end{table}

\subsection{Two-dimensional density current}
The computational domain in the $xz$-plane is $\Omega = [0,25600] \times [0,6400] \, \mathrm{m^2}$ and the time interval of interest is $(0,  900 ]$ s. Impenetrable, free-slip boundary conditions are imposed on the whole domain boundary. The initial potential temperature,  $\theta^0$, is
\begin{equation}
\theta^0 = \left\{
\begin{array}{ll}
 300 - \displaystyle \frac{15}{2} \left[ 1 + \cos(\pi r) \right] \qquad &\text{if } r \leq 1 \, \mathrm{m}\\[2mm] 300 &\text{ otherwise},
\end{array}
\right.
\end{equation}
with $r = \sqrt{\left(\frac{x - x_c}{x_r}\right)^2 + \left(\frac{z - z_c}{z_r}\right)^2}$ for $(x_r, z_r) = (4000, 2000) \, \mathrm{m}$ and $(x_c, z_c) = (0, 3000) \, \mathrm{m}$.
In this benchmark, the initial perturbation is a bubble of cold air, which descends due to negative buoyancy and, when it reaches the ground, it rolls up and forms a cold front. The initial density $\rho^0$ is computed by solving the hydrostatic balance (\ref{eq:14}), while the initial velocity field is zero everywhere. Moreover, in accordance with \cite{Restelli2009}, 
we set $\mu_a = 75$ in \eqref{eq:21} and $Pr = 1$ in \eqref{eq:22}.

We consider three meshes with a uniform resolution $h=\Delta x = \Delta z=$ 50, 100, 200 m, an a time step $\Delta t =0.1$ s. Following the approach adopted for the previous test, we use as reference the solution computed on the finest grid without mesh adaptation, consistently with the extensive validation of the underlying numerical model in \cite{GIRFOGLIO2025106510,GQR_OF_clima,CGQR,Clinco2024}.

For a fair comparison between the IREE and PMA algorithms, we set
$I_r=5$ and $N_\text{max}=2\times10^6$ for both approaches.
All the other parameters are set according to the values in Table~\ref{tab:9}. 
\begin{table}[ht]
    \centering
    
        \begin{minipage}{0.45\textwidth}  
        \centering
        \begin{tabular}{|c|c|c|c|c|}
            \hline
            \multicolumn{5}{|c|}{\textbf{IREE}} \\ \hline
            \textbf{Mesh size} & $\delta_1$ & $\delta_2$ & $tol$ & $l_\text{max}$\\ \hline
            100 & 8.0 & 0.99 & 1e-1 & 1  \\ \hline
            200 & 6.0 & 0.99 & 1e-1 & 2  \\ \hline
        \end{tabular}
    \end{minipage}
    \begin{minipage}{0.45\textwidth}  
        \centering
        \begin{tabular}{|c|c|c|c|c|}
            \hline
            \multicolumn{5}{|c|}{\textbf{PMA}} \\ \hline
            \textbf{Mesh size} & $\alpha_\text{min}^\text{ref}$ & $\alpha_\text{max}^\text{ref}$ & $l_\text{max}$ & $l_c$ \\ \hline
            100 & 0.1 & 0.9 & 1 & 10  \\ \hline
            200 & 0.1 & 0.9 & 2 & 10  \\ \hline
        \end{tabular}
    \end{minipage}%
                \caption{2D density current: parameter values for the IREE and PMA algorithms.}\label{tab:9}
\end{table}

\begin{figure}[h!]
    \centering
    
    \begin{overpic}[width=0.48\textwidth]{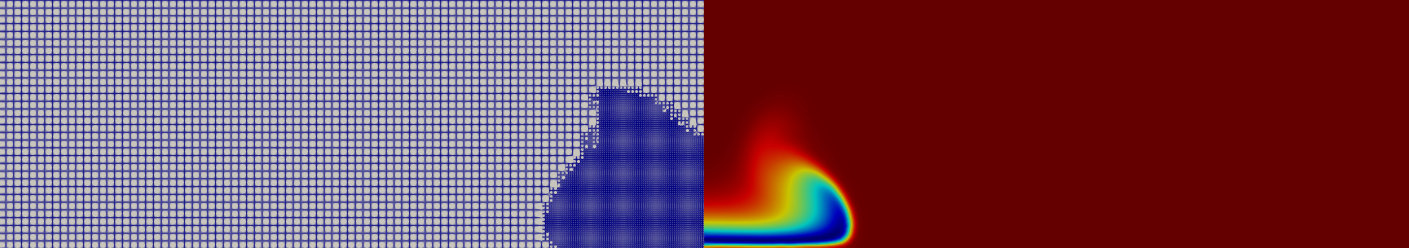}
        \put(70,14){\textcolor{white}{\footnotesize{$t = 300$ s, IREE}}}
    \end{overpic}
    \begin{overpic}[width=0.48\textwidth]{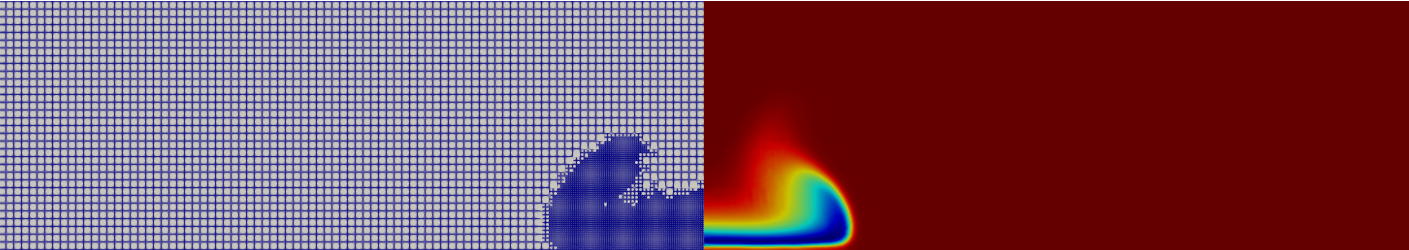}
        \put(70,14){\textcolor{white}{\footnotesize{$t = 300$ s, PMA}}}
    \end{overpic}
    \\
    \vspace{3pt}
    \begin{overpic}[width=0.48\textwidth]{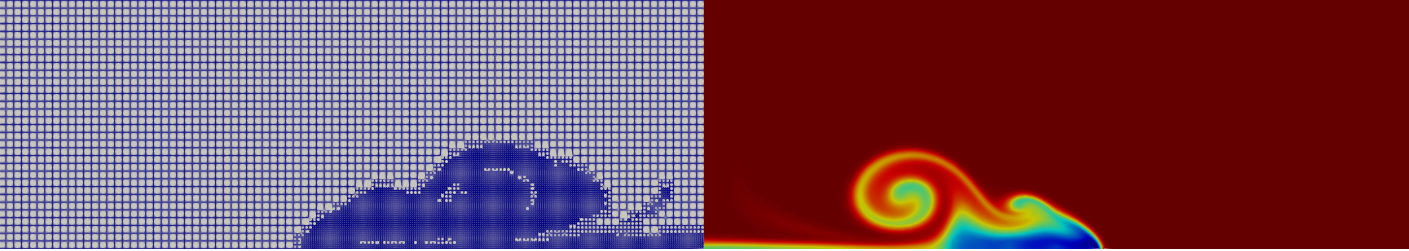}
        \put(70,14){\textcolor{white}{\footnotesize{$t = 600$ s, IREE}}}
    \end{overpic}
    \begin{overpic}[width=0.48\textwidth]{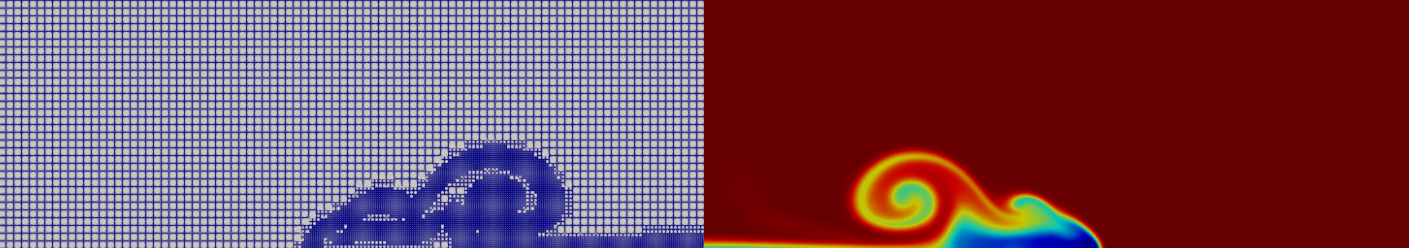}
        \put(70,14){\textcolor{white}{\footnotesize{$t = 600$ s, PMA}}}
    \end{overpic}
    \\
    \vspace{3pt}
    
    \begin{overpic}[width=0.48\textwidth]{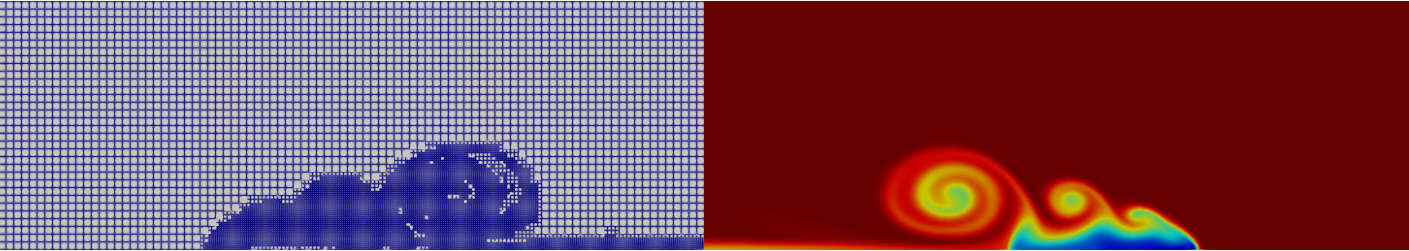}
        \put(70,14){\textcolor{white}{\footnotesize{$t = 750$ s, IREE}}}
    \end{overpic}
    \begin{overpic}[width=0.48\textwidth]{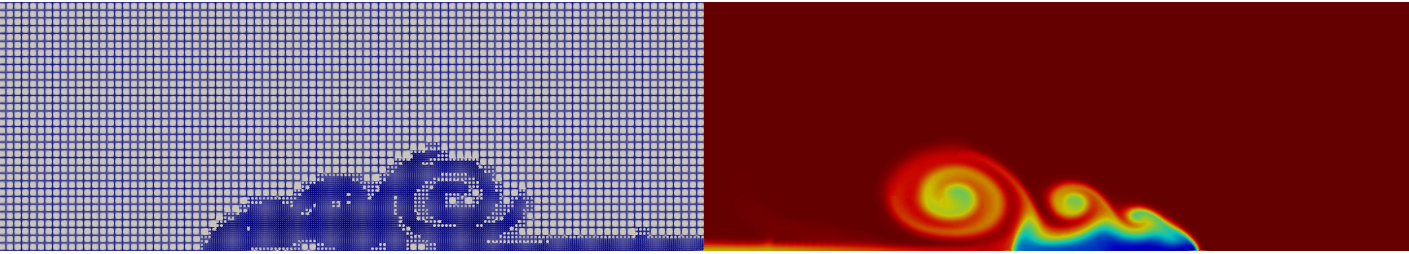}
        \put(70,14){\textcolor{white}{\footnotesize{$t = 750$ s, PMA}}}
    \end{overpic}
    \\
    \vspace{3pt}
    
    \begin{overpic}[width=0.48\textwidth]{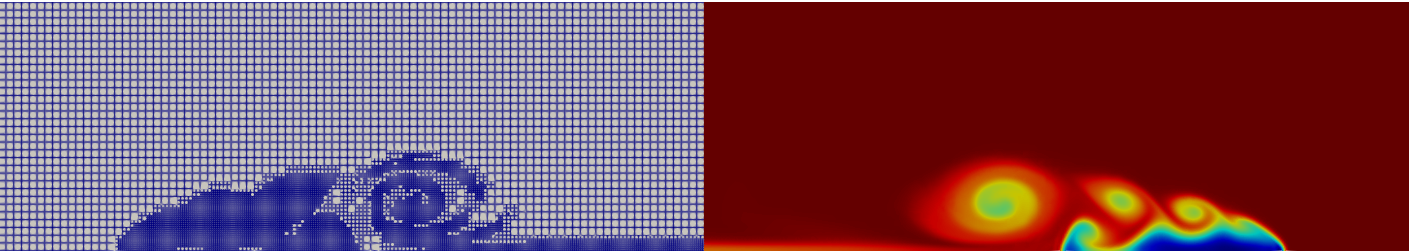}
        \put(70,14){\textcolor{white}{\footnotesize{$t = 900$ s, IREE}}}
    \end{overpic}
    \begin{overpic}[width=0.48\textwidth]{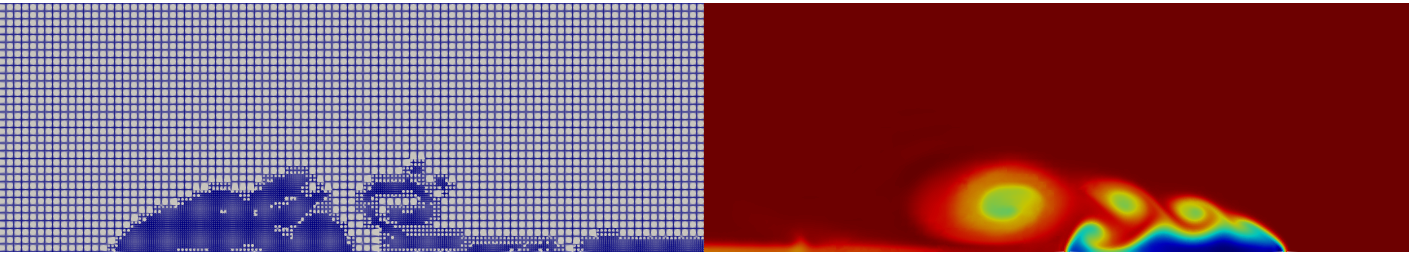}
        \put(70,14){\textcolor{white}{\footnotesize{$t = 900$ s, PMA}}}
    \end{overpic}

    \vspace{7pt}

            \begin{overpic}[width=0.60\textwidth]{Images/colorbar/colorbar.png}
        \put(53,85){\textcolor{white}{\footnotesize{}}}
    \end{overpic}
    \caption{2D density current: 
    evolution of the potential temperature and associated mesh adapted by the IREE (left) and PMA (right) algorithms starting from a mesh size $h=200$ m. }
    \label{fig:DensityCurrent2DMesh}
\end{figure}

Figure~\ref{fig:DensityCurrent2DMesh} shows the evolution of the computed potential temperature and the corresponding
mesh adapted by the IREE and PMA algorithms,
starting from a mesh size $h=200$ m. 
As in the previous test, the IREE algorithm refines
a larger portion of the computational domain than the PMA algorithm, although the difference between the two strategies progressively decreases as the simulation evolves.

\begin{figure}[h!]
    \centering
    
    \begin{overpic}[width=0.32\textwidth]{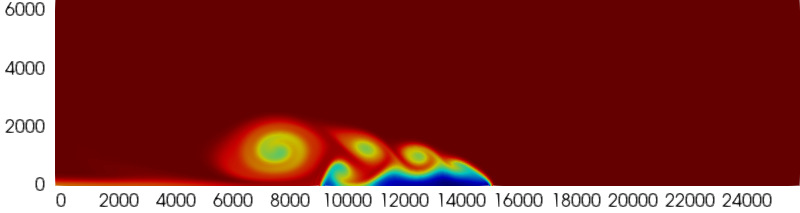}
        \put(25,20){\textcolor{white}{\footnotesize{Reference ($h = 50$ m)}}}
    \end{overpic}
    \begin{overpic}[width=0.32\textwidth]{Images/DC_contour/T50_contour.png}
        \put(25,20){\textcolor{white}{\footnotesize{Reference ($h = 50$ m)}}}
    \end{overpic}
    \begin{overpic}[width=0.32\textwidth]{Images/DC_contour/T50_contour.png}
        \put(25,20){\textcolor{white}{\footnotesize{Reference ($h = 50$ m)}}}
    \end{overpic}
    \\
    
    \begin{overpic}[width=0.32\textwidth]{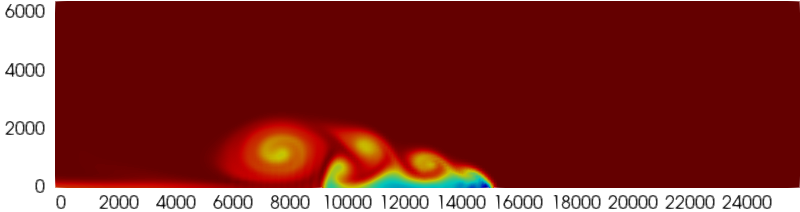}
        \put(24,20){\textcolor{white}{\footnotesize{No Adapt., $h = 100$ m}}}
    \end{overpic}
        \begin{overpic}[width=0.315\textwidth]{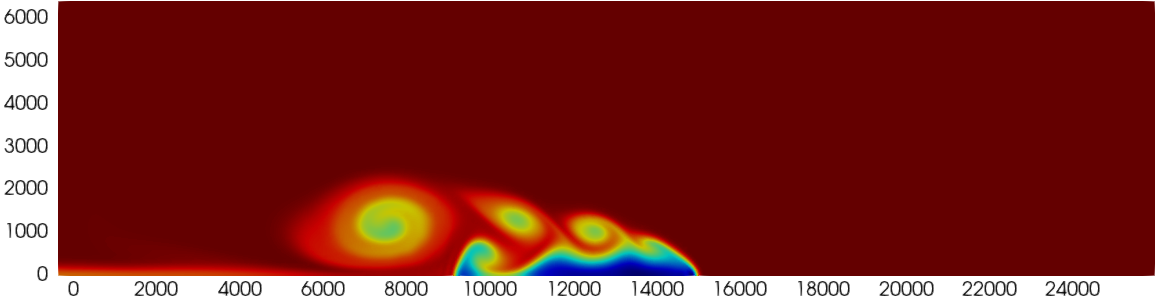}
        \put(30,20){\textcolor{white}{\footnotesize{IREE, $h = 100$ m}}}
    \end{overpic}
    \begin{overpic}[width=0.32\textwidth]{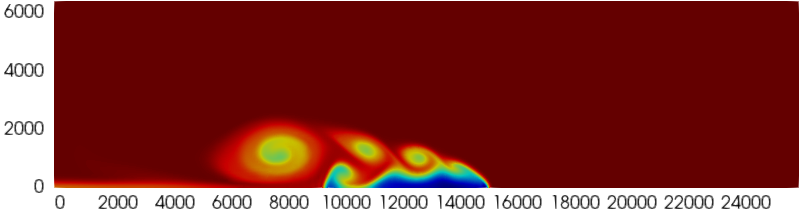}
        \put(30,20){\textcolor{white}{\footnotesize{PMA, $h = 100$ m}}}
    \end{overpic}
    \\
    
    \begin{overpic}[width=0.32\textwidth]{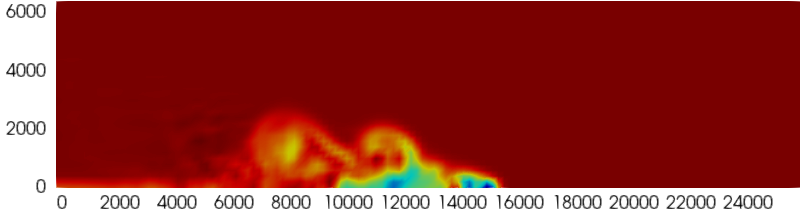}
        \put(24,20){\textcolor{white}{\footnotesize{No Adapt., $h = 200$ m}}}
    \end{overpic}
        \begin{overpic}[width=0.32\textwidth]{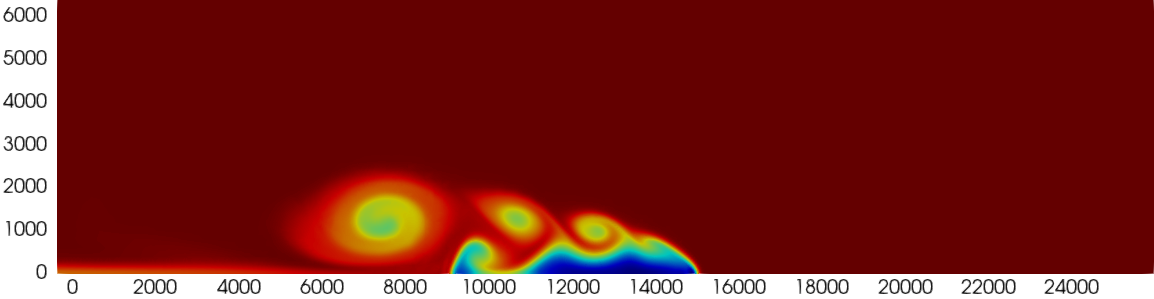}
        \put(30,20){\textcolor{white}{\footnotesize{IREE, $h = 200$ m}}}
    \end{overpic}
    \begin{overpic}[width=0.32\textwidth]{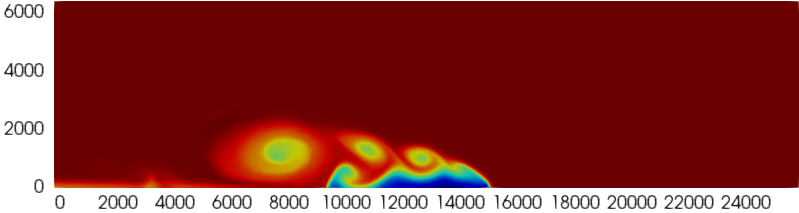}
        \put(30,20){\textcolor{white}{\footnotesize{PMA, $h = 200$ m}}}
    \end{overpic}
    \\
    \vspace{7pt}

            \begin{overpic}[width=0.50\textwidth]{Images/colorbar/colorbar.png}
    \end{overpic}
    
    \caption{2D density current:  potential temperature at time $t=900$ s provided by the reference simulation (first row) compared against the solution computed with mesh size $h = 100$ m (second row), and $h = 200$ m (third row), with no mesh adaptation (first column) and with mesh adaptation driven by the IREE (second column) and PMA (third column) algorithm.
    }
    \label{fig:DensityCurrent2Dcontour}
\end{figure}

Figure~\ref{fig:DensityCurrent2Dcontour} compares the potential temperature computed at $t=900$ s, at the end of the selected time interval, for all the meshes under consideration, with and without mesh adaptation.
As expected, the solution shown in the first column and computed without mesh adaptation 
deteriorates drastically in accuracy as the mesh becomes coarser.
The mesh adaptation algorithms allow for a substantial recovery of accuracy for both the choices $h = 100, 200$ m, by capturing the tri-rotor structure that emerges at $t=900$ s. 
We observe that, for $h=200$ m, the IREE algorithm yields a smoother solution in the wake of the cold front, whereas the PMA solution still exhibits localized instabilities.

For a quantitative comparison, in Fig.~\ref{fig:DensityCurrent2Dcontour} we provide the evolution of the $L^2(\Omega)$-norm of the relative error for all the considered mesh resolutions and adaptation strategies. As in Fig.~\ref{fig:RTB_error}, the simulations without mesh adaptation yield
the largest errors. The adoption of mesh adaptation improves the solution accuracy by approximately one order of magnitude. For both mesh sizes, the PMA algorithm initially provides the most accurate results, achieving the lowest errors up to approximately $t=300$ s. Beyond this stage, the IREE algorithm becomes progressively more competitive and eventually yields the smallest errors for the majority of the simulation. This suggests that PMA is preferable when short-time accuracy is the primary objective, whereas IREE offers superior performance in long-time simulations.
\begin{figure}
    \centering
    \begin{subfigure}{0.48\textwidth}
        \includegraphics[width = \textwidth]{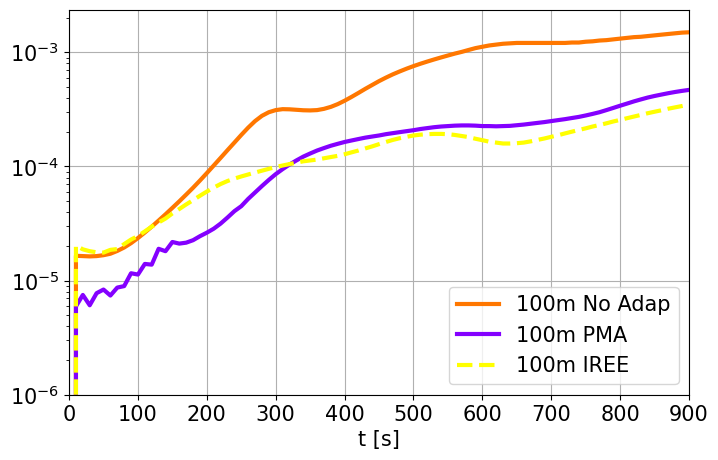}
    \end{subfigure}
    \begin{subfigure}{0.48\textwidth}
        \includegraphics[width = \textwidth]{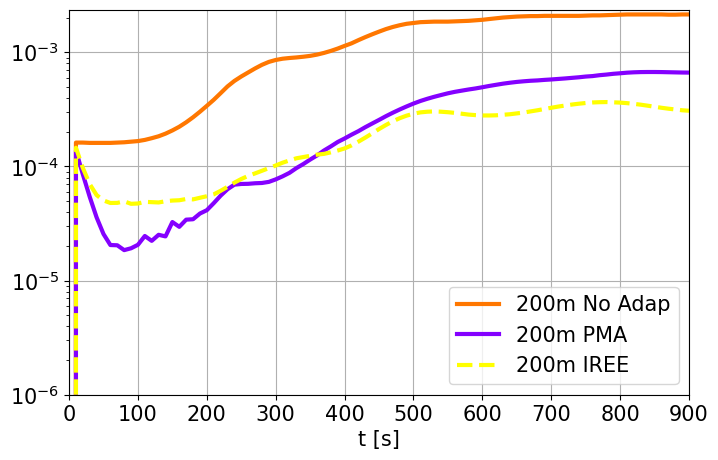}
    \end{subfigure}
    \caption{2D density current: the time evolution of the relative error in the $L^2(\Omega)$-norm  for the potential temperature $\theta$ using mesh sizes $h = 100$ m (left) and $h = 200$ m (right). The reference solution in the calculation of the error is the approximation computed on a fixed uniformly fine grid with mesh size $h = 50$ m.}
    \label{fig:DC2D_error}
\end{figure}

\begin{figure}[h!]
    \centering
    \begin{subfigure}{0.48\textwidth}
        \centering
        \includegraphics[width=\textwidth]{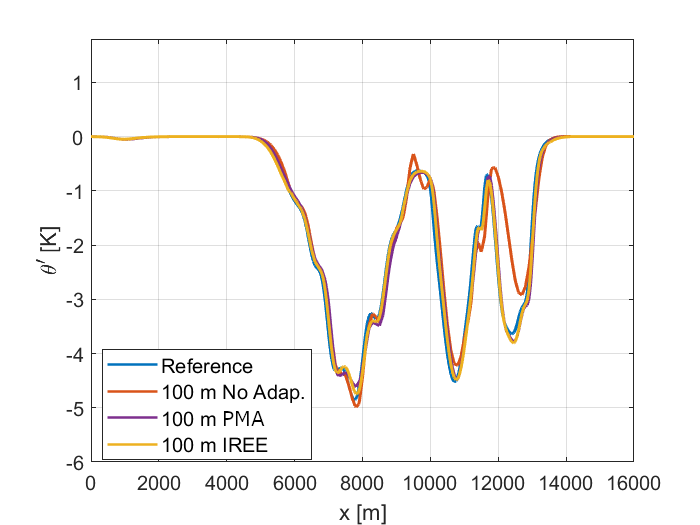}
        \label{fig:immagine1c}
    \end{subfigure}
    \begin{subfigure}{0.48\textwidth}
        \centering
        \includegraphics[width=\textwidth]{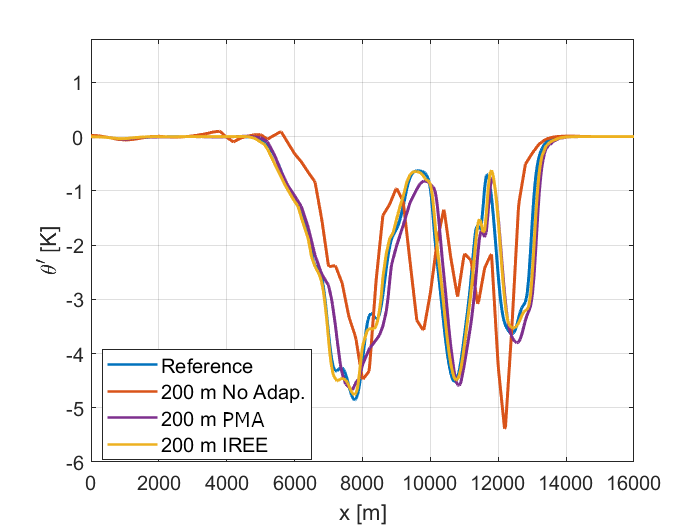}
        \label{fig:immagine2c}
    \end{subfigure}
    \caption{2D density current: potential temperature flluctuation, $\theta'$, at $t = 900$ s
along the horizontal direction at a height of $z = 1200$ m, computed with a mesh size 
$h = 100$ m (left) and $h = 200$ m (right). The reference values are associated with a fixed uniformly fine grid with mesh size $h = 50$ m.}
    \label{fig:DensityCurrent2Dmatlab}
\end{figure}

A different perspective on the comparison between IREE and PMA is provided by Fig.~\ref{fig:DensityCurrent2Dmatlab} which shows the fluctuations,  $\theta'$, of the potential temperature along the horizontal $x$-direction at a fixed height $z=1200$ m at time $t=900$ s, following the analysis proposed in  \cite{Giraldo2008}. 
In both plots of Fig.~\ref{fig:DensityCurrent2Dmatlab}, the solution obtained without mesh adaptation progressively departs from the reference solution, with the discrepancy increasing over time and becoming more pronounced on the coarser mesh. In contrast, the adaptive simulations remain much closer to the reference curve, with the IREE algorithm consistently providing the best agreement. For the finer mesh ($h=100$ m), both adaptive strategies reproduce the reference dynamics reasonably well. For the coarser mesh ($h=200$ m), however, noticeable differences emerge: while the IREE solution remains in good agreement with the reference evolution, the PMA solution progressively loses accuracy, exhibiting both amplitude attenuation and a slight phase shift in the prediction of successive extrema.

Finally, we investigate the computational efficiency of the two adaptive strategies. The reference simulation requires $734$ s, whereas the computational times associated with meshes $h=100,200$ m and the IREE and PMA algorithms are reported in Table~\ref{tab:10}, together with the corresponding percentage time savings. 
The results confirm that both adaptive approaches lead to substantial reductions in computational cost with respect to the reference simulation. As expected, PMA remains the less expensive strategy. However, unlike the rising thermal bubble test, the difference in computational time between IREE and PMA is less pronounced. This behavior is consistent with the refinement patterns shown in Fig.~\ref{fig:DensityCurrent2DMesh}: although IREE still refines a larger portion of the computational domain, the discrepancy between the regions refined by the two algorithms is smaller than in the previous benchmark. 
The moderate increase in computational cost associated with IREE is compensated by its improved accuracy, particularly for the coarser mesh, as highlighted by the error trends discussed above and by the results reported in Fig.~\ref{fig:DC2D_error}.
\begin{table}[h!]
\centering
\begin{tabular}{|c||c|c|}
\hline
mesh size $h$ [m]& IREE algorithm [s]& PMA algorithm [s] \\
\hline
100 &  471 (35.8\%) & 382 (47.9\%)\\                 
200 & 231 (68.5\%)  & 153 (79.1\%) \\                 
\hline
\end{tabular}
\caption{2D density current: 
computational performance of the IREE and PMA algorithms for mesh sizes $h = 100, 200$ m, in terms of computational time and percentage of time saved with respect to the reference simulation ($734$ s).
}\label{tab:10}
\end{table}

\subsection{Three-dimensional density current}
The computational domain for this benchmark is $\Omega = [0,12000] \times [0,6400] \times [0, 12000] \, \mathrm{m^3}$ and the time interval of interest is $(0, 600]$ s. Impenetrable, free-slip boundary conditions are imposed on the domain walls. Like in the two-dimensional case, we start from a neutrally stratified atmosphere with uniform background potential temperature $\theta^0 = 300\,\,K$, perturbed by a spherical bubble of colder air, being 
\begin{equation}
\theta^0 = \left\{
\begin{array}{ll}
 300 - \displaystyle \frac{15}{2} \left[ 1 + \cos(\pi r) \right] \qquad &\text{if } r \leq 1 \, \mathrm{m}\\[2mm] 300 &\text{ otherwise},
\end{array}
\right.
\end{equation}
with $r = \sqrt{\left(\frac{x - x_c}{x_r}\right)^2 +\left(\frac{y - y_c}{y_r}\right)^2+ \left(\frac{z - z_c}{z_r}\right)^2}$ for $(x_r, y_r, z_r) = (4000,2000, 4000) \, \mathrm{m}$ and $(x_c, y_c, z_c) = (0, 3000, 0) \, \mathrm{m}$. The initial density $\rho^0$ is computed by solving the hydrostatic balance (\ref{eq:14}), while the initial velocity field is zero everywhere. For the artificial viscosity and the Prandtl number, we use the same values as in the 2D cases, by setting $\mu_a = 75$, and $Pr = 1$.

We consider three meshes with uniform resolution $h=\Delta x = \Delta y= \Delta z=$ 50, 200, 400 m, and a time step $\Delta t=0.1$ s. As in the previous test cases, we adopt as reference
the discrete solution computed on the finest grid without mesh adaptation.
The underlying mesh contains approximately 6.9 million cells, making the corresponding simulation particularly demanding from a computational standpoint and resulting in a runtime of about 58 hours. 

For the comparison between the IREE and PMA algorithms, the common parameters
$I_r = 5$ and $N_\text{max}=2\times10^8$ are used in both approaches. The remaining parameter values are reported in Table~\ref{tab:12}.
%
\begin{table}[ht]
    \centering 
        \begin{minipage}{0.45\textwidth}  
        \centering
        \begin{tabular}{|c|c|c|c|c|}
            \hline
            \multicolumn{5}{|c|}{\textbf{IREE}} \\ \hline
            \textbf{Mesh size} & $\delta_1$ & $\delta_2$ & $tol$ & $l_\text{max}$\\ \hline
            200 & 6.5 & 0.99 & 1e-0 & 2  \\ \hline
            400 & 2.1 & 0.99 & 1e-0 & 3  \\ \hline
        \end{tabular}
    \end{minipage}
    \begin{minipage}{0.45\textwidth}  
        \centering
        \begin{tabular}{|c|c|c|c|c|}
            \hline
            \multicolumn{5}{|c|}{\textbf{PMA}} \\ \hline
            \textbf{Mesh size} & $\alpha_\text{min}^\text{ref}$ & $\alpha_\text{max}^\text{ref}$ & $l_\text{max}$ & $l_c$ \\ \hline
            200 & 0.1 & 0.9 & 2 & 10  \\ \hline
            400 & 0.1 & 0.9 & 3 & 10  \\ \hline
        \end{tabular}
    \end{minipage}%
        \caption{3D density current: parameter values for the IREE and PMA algorithms.}\label{tab:12}
\end{table}

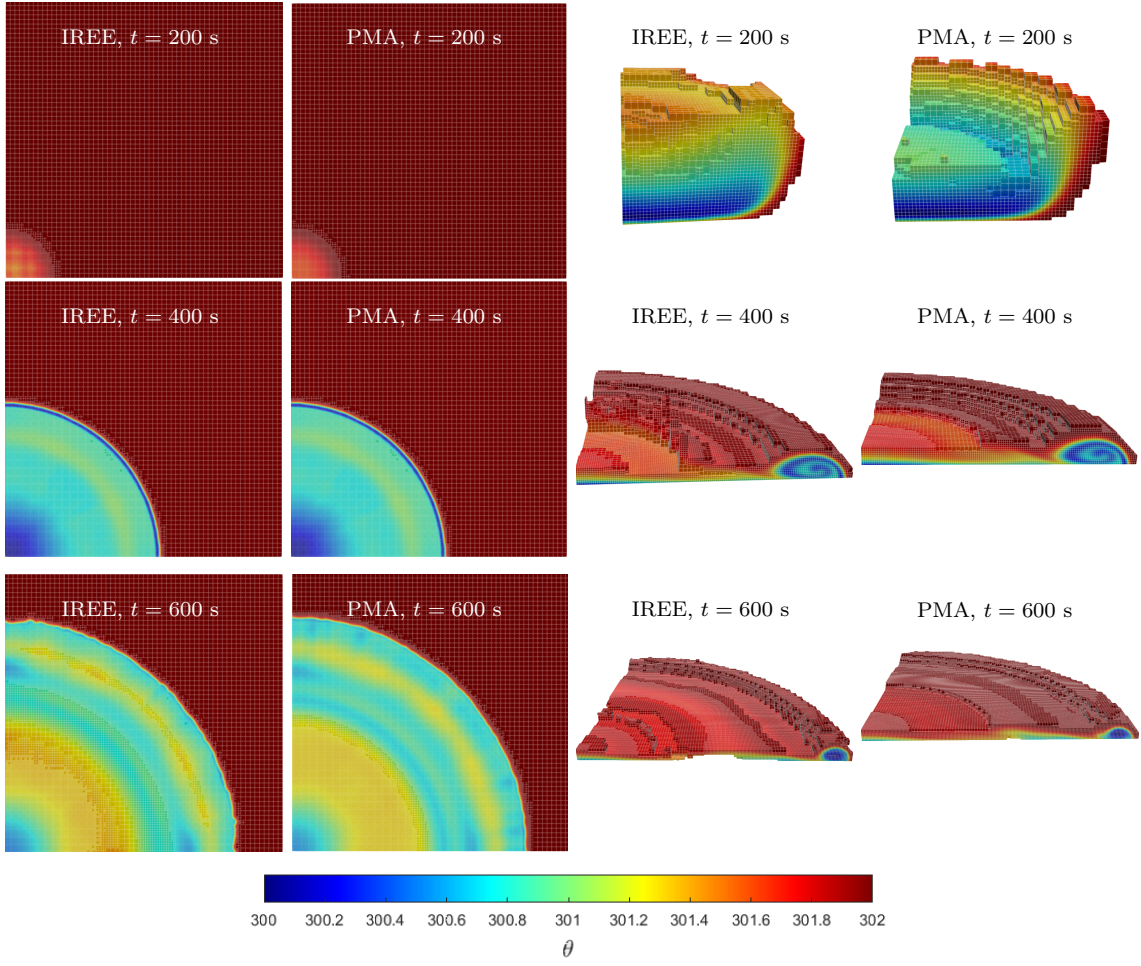
\begin{figure}[h!]
\centering
    \begin{minipage}{0.24\textwidth}
        \begin{overpic}[width=\linewidth]{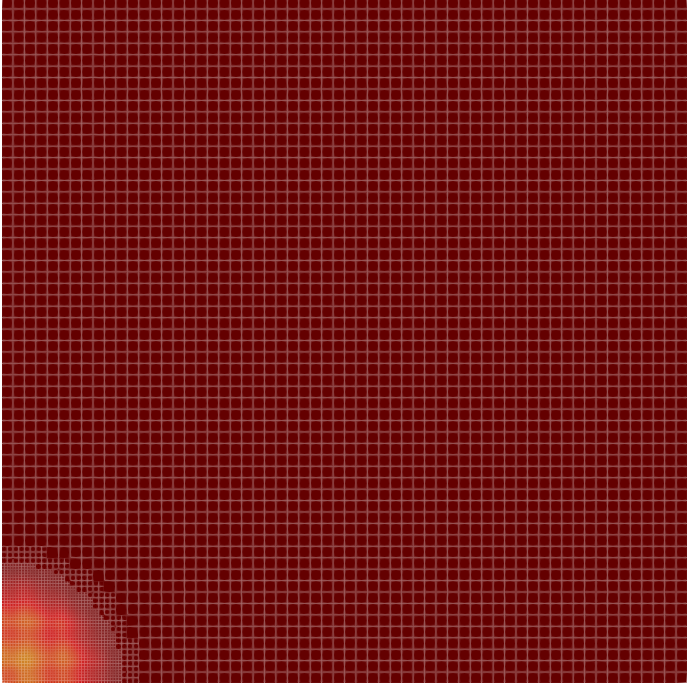}
            \put(20,85){\textcolor{white}{\footnotesize{IREE, $t = 200$ s}}}
            \put(65,103){\textcolor{black}{{}}}
        \end{overpic}
    \end{minipage}
    \begin{minipage}{0.24\textwidth}
        \begin{overpic}[width=\linewidth]{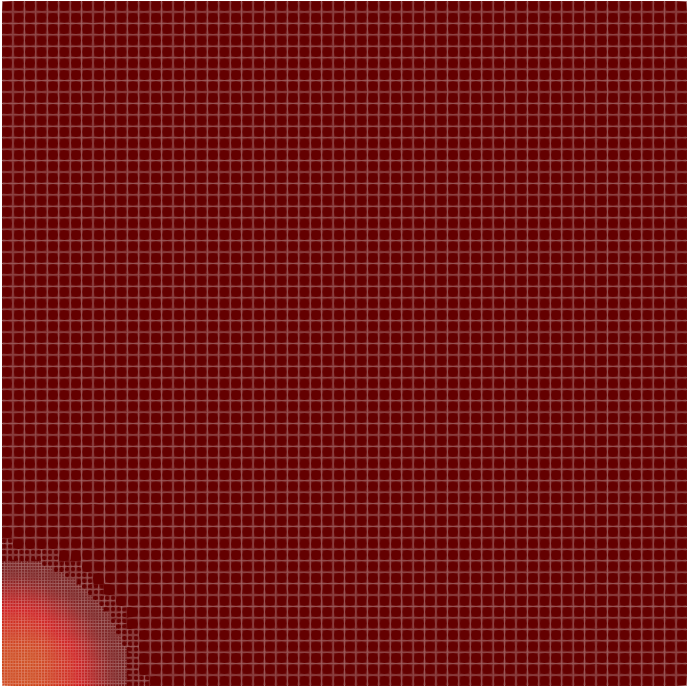}
            \put(20,85){\textcolor{white}{\footnotesize{PMA, $t = 200$ s}}}
        \end{overpic}
    \end{minipage}
    \begin{minipage}{0.24\textwidth}
        \begin{overpic}[width=\linewidth]{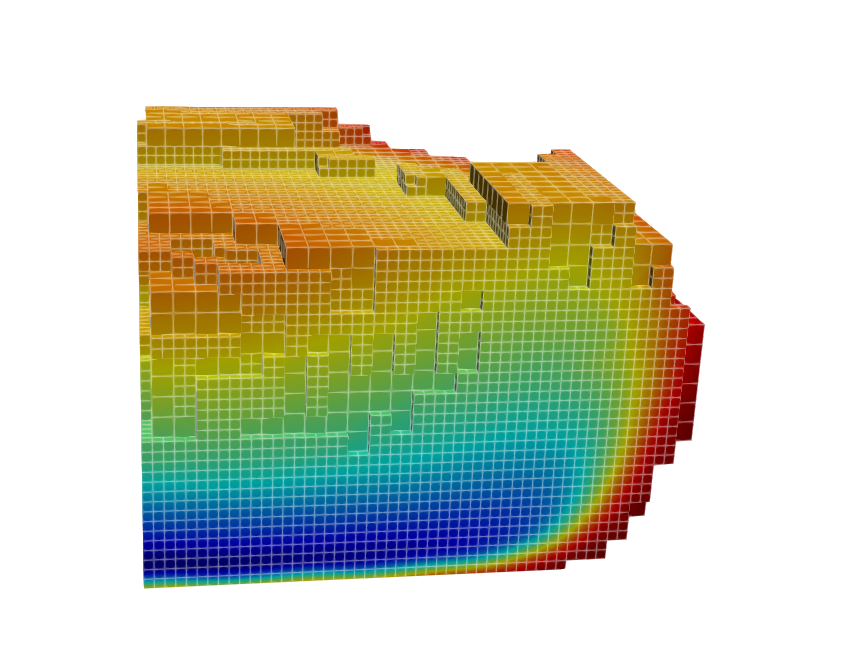}
            \put(88,90){\textcolor{black}{{}}}
            \put(20,74){\textcolor{black}{\footnotesize{IREE, $t = 200$ s}}}
        \end{overpic}
    \end{minipage}
    \begin{minipage}{0.24\textwidth}
        \begin{overpic}[width=\linewidth]{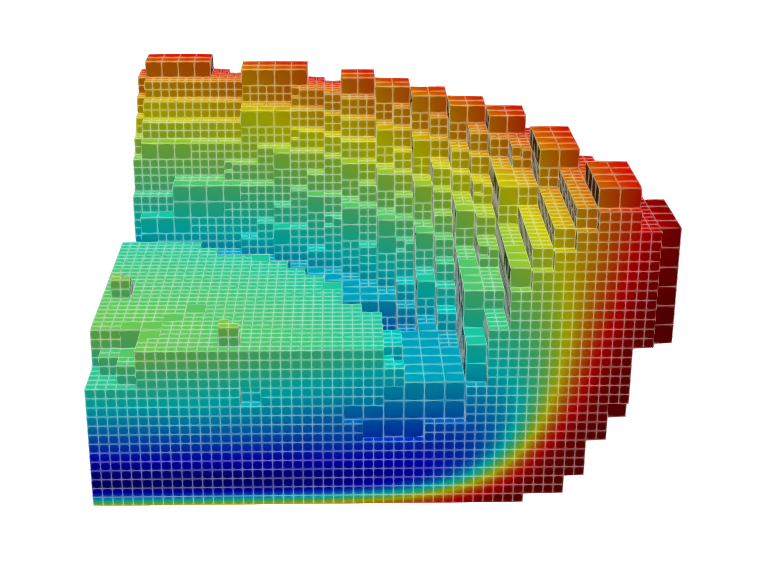}
\put(20,73){\textcolor{black}{\footnotesize{PMA, $t = 200$ s}}}
        \end{overpic}
    \end{minipage}
    \vspace{5pt} 
    \begin{minipage}{0.24\textwidth}
        \begin{overpic}[width=\linewidth]{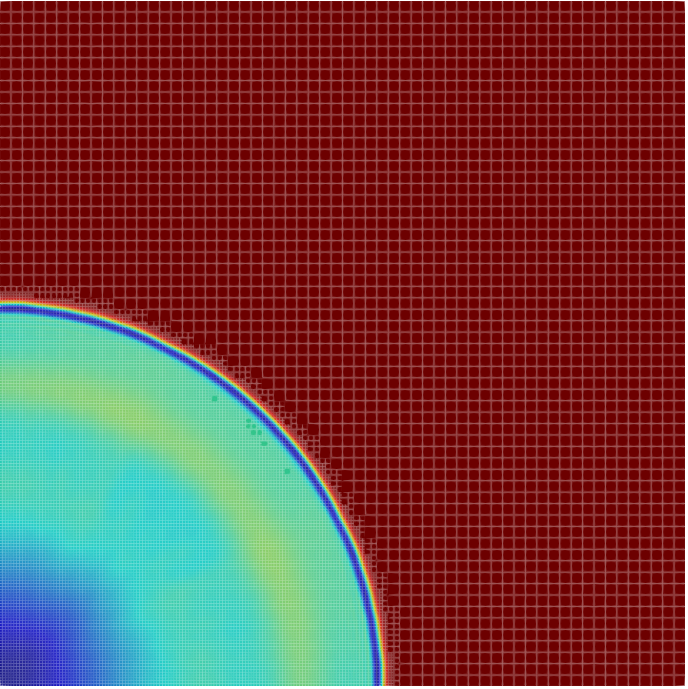}
          \put(20,85){\textcolor{white}{\footnotesize{IREE, $t = 400$ s}}}
        \end{overpic}
    \end{minipage}
    \begin{minipage}{0.24\textwidth}
        \begin{overpic}[width=\linewidth]{Images/DC3D/DC3D_impronta_IREE400.png}
        \put(20,85){\textcolor{white}{\footnotesize{PMA, $t = 400$ s}}}
        \end{overpic}
    \end{minipage}
    \begin{minipage}{0.24\textwidth}
        \begin{overpic}[width=\linewidth]{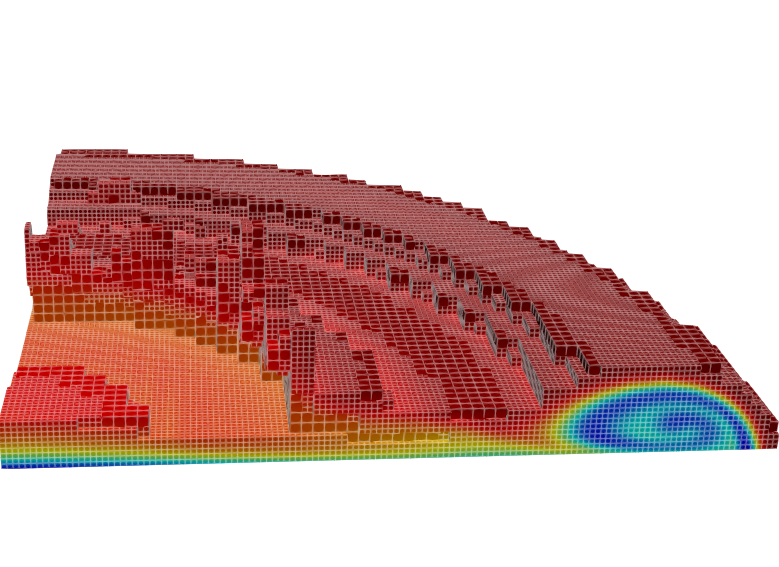}
        \put(20,71){\textcolor{black}{\footnotesize{IREE, $t = 400$ s}}}
        \end{overpic}
    \end{minipage}
    \begin{minipage}{0.24\textwidth}
        \begin{overpic}[width=\linewidth]{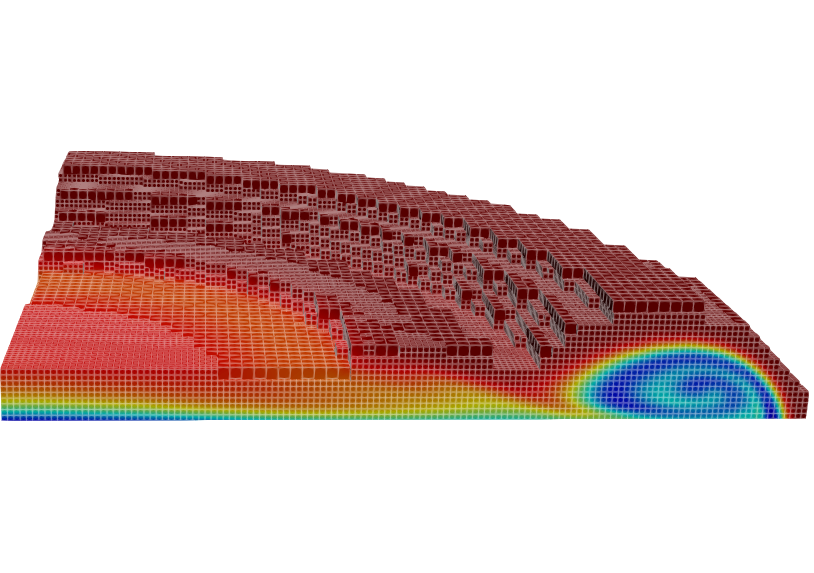}
        \put(20,70){\textcolor{black}{\footnotesize{PMA, $t = 400$ s}}}
        \end{overpic}
    \end{minipage}
    \vspace{5pt} 
    \begin{minipage}{0.24\textwidth}
        \begin{overpic}[width=\linewidth]{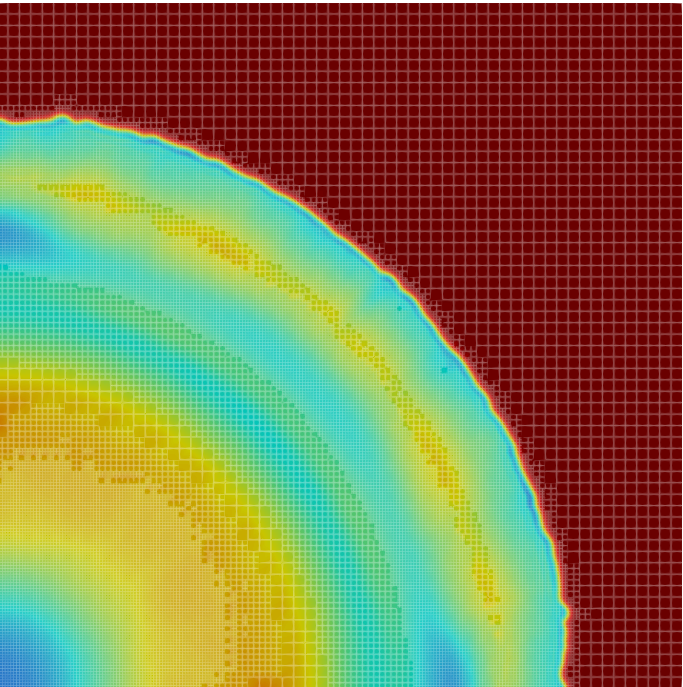}
            \put(20,85){\textcolor{white}{\footnotesize{IREE, $t = 600$ s}}}
        \end{overpic}
    \end{minipage}
    \begin{minipage}{0.24\textwidth}
        \begin{overpic}[width=\linewidth]{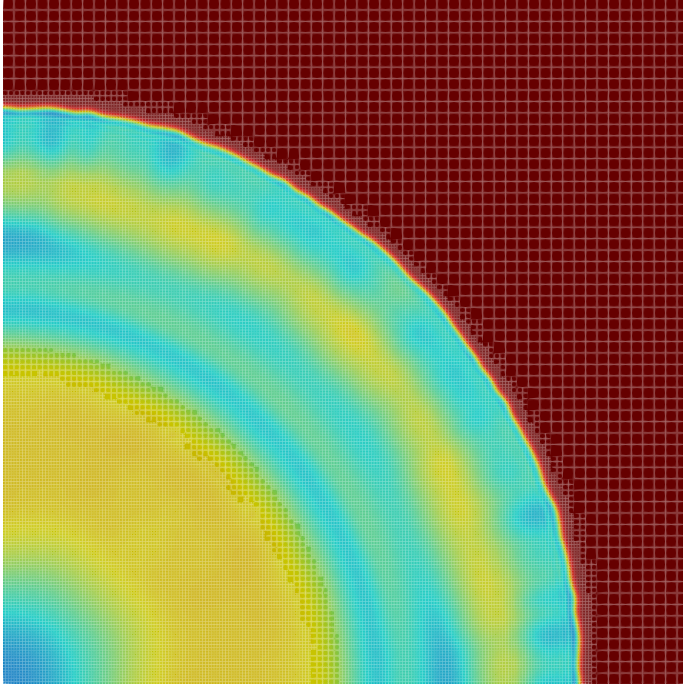}
        \put(20,85){\textcolor{white}{\footnotesize{PMA, $t = 600$ s}}}
        \end{overpic}
    \end{minipage}
    \begin{minipage}{0.24\textwidth}
        \begin{overpic}[width=\linewidth]{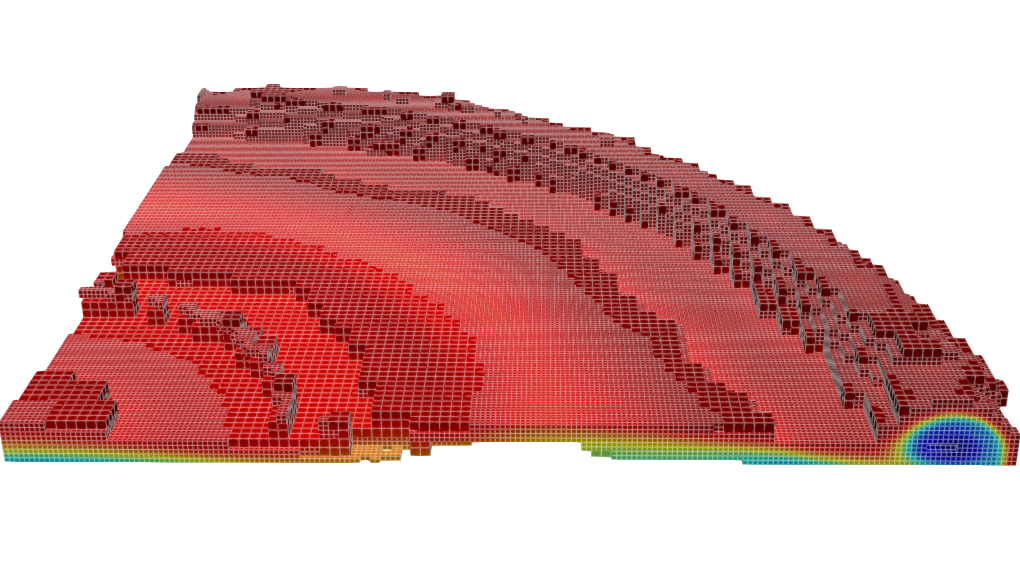}
        \put(20,63){\textcolor{black}{\footnotesize{IREE, $t = 600$ s}}}
        \end{overpic}
    \end{minipage}
    \begin{minipage}{0.24\textwidth}
        \begin{overpic}[width=\linewidth]{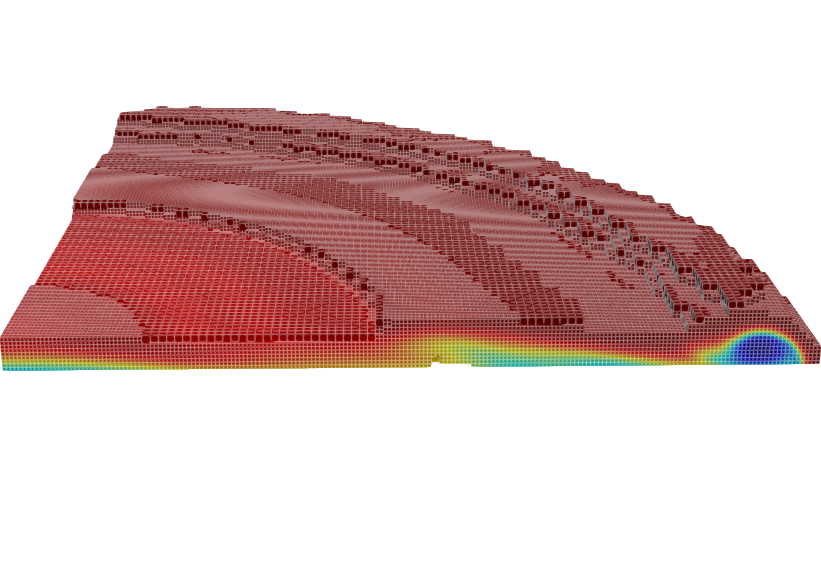}
        \put(20,69){\textcolor{black}{\footnotesize{PMA, $t = 600$ s}}}
        \end{overpic}
    \end{minipage}
    \vspace{8pt}
            \begin{overpic}[width=0.60\textwidth]{Images/colorbar/colorbar.png}
    \end{overpic}
    \caption{3D density current: 
    evolution of the potential and associated mesh adapted by the IREE (first and third column) and the PMA (second and fourth column) algorithms starting from a mesh size $h=200$ m. 
The first two columns show the solution at ground level ($z=0$), whereas the last two columns display a three-dimensional cross-sectional view of the domain.}
    \label{fig:DC3DInTime}
\end{figure}

Figure~\ref{fig:DC3DInTime} shows the evolution of
the computed potential temperature 
and details of the meshes 
given by the IREE and PMA algorithms. Specifically, 
the first two columns display the potential temperature at the ground level (i.e., for $z = 0$), while the last two columns show the set of the cells that have been refined. As in the 2D benchmarks, the
IREE algorithm refines a larger portion of the computational domain than the PMA algorithm. 

\begin{figure}[htb!]
    \centering
    \begin{subfigure}{0.48\textwidth}
        \includegraphics[width = \textwidth]{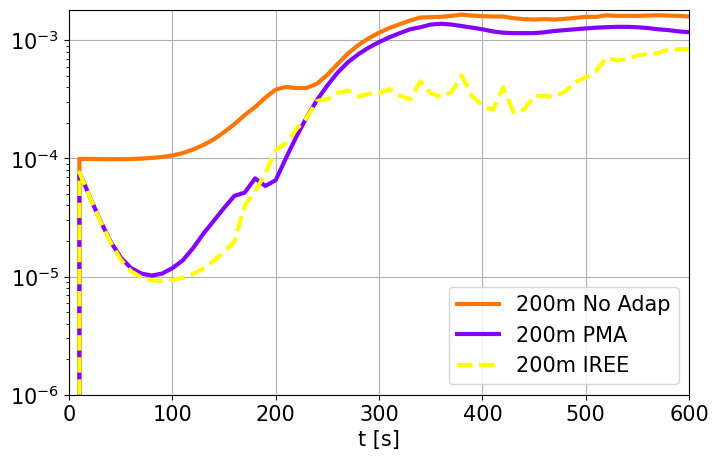}
    \end{subfigure}
    \begin{subfigure}{0.48\textwidth}
        \includegraphics[width = \textwidth]{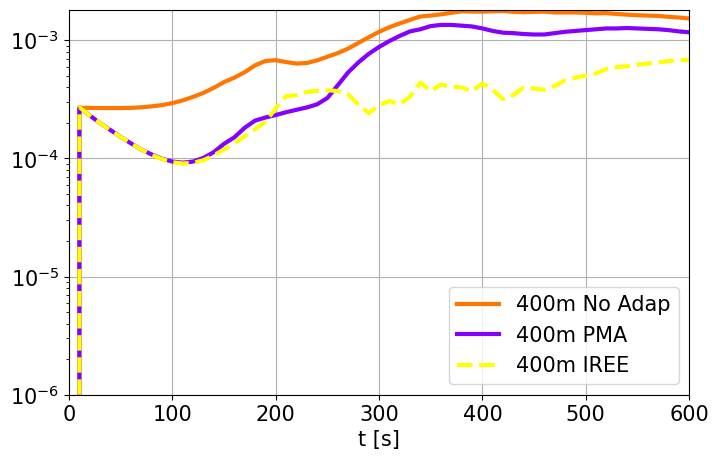}
    \end{subfigure}
    \caption{3D density current: the time evolution of the relative error in the $L^2(\Omega)$-norm for the potential temperature $\theta$ using mesh sizes $h = 200$ m (left) and $h = 400$ m (right). The reference solution in the calculation of the error is the approximation computed on a fixed uniformly fine grid with $h = 50$ m.}
    \label{fig:DC3D_error}
\end{figure}

The next step is to assess whether this additional refinement results in a more accurate solution. Figure~\ref{fig:DC3D_error} shows the evolution of the $L^2(\Omega)$-norm of the relative error characterizing the simulations with mesh sizes $h = 200, 400$ m, with and without mesh adaptation. As we observed in Fig.~\ref{fig:RTB_error} and Fig.~\ref{fig:DC2D_error}, the simulations with coarse meshes and no mesh adaptation give the largest
errors. Mesh adaptation algorithms reduce the error by up to an order of magnitude. For both meshes $h = 200, 400$ m, the PMA algorithm slightly outperforms the IREE algorithm only over the time interval $[200,250]$, whereas the IREE algorithm achieves higher accuracy throughout the remainder of the simulation (approximately 92\% of the simulated time).

\begin{figure}
    \centering
    \begin{subfigure}{0.48\textwidth}
        \centering
        \includegraphics[width=\textwidth]{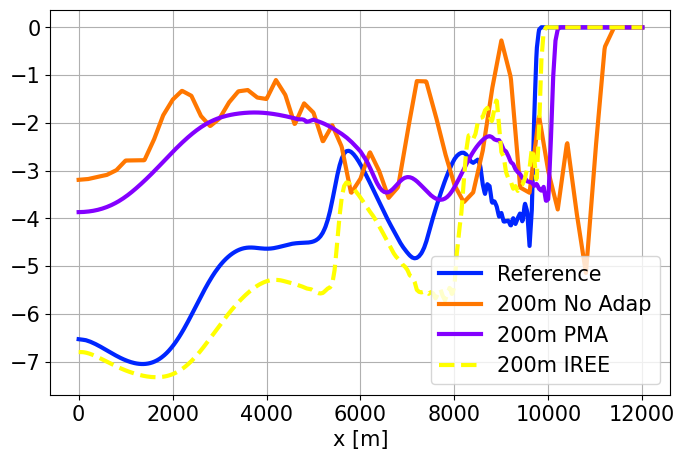}
        \label{fig:immagine1d}
    \end{subfigure}
    \begin{subfigure}{0.48\textwidth}
        \centering
        \includegraphics[width=\textwidth]{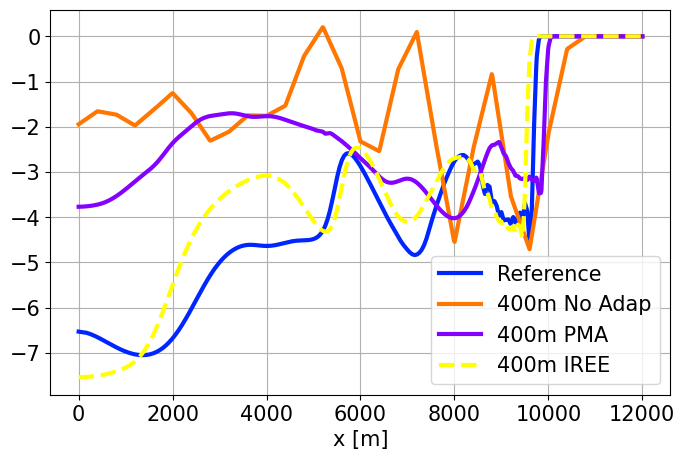}
        \label{fig:immagine2e}
    \end{subfigure}
    \caption{3D density current: potential temperature fluctuation, $\theta'$, at $t = 600$ s at a height of $z = 280$ m, along the line connecting points $(0,0,280)$ to $(12000,0,280)$, computed with mesh size $h = 200$ m (left), $h = 400$ m (right). The reference values are associated with a fixed uniformly fine grid with $h = 50$ m.}
    \label{fig:DC3D_matlab} 
\end{figure}

To further investigate the behavior of the two adaptive strategies, in Fig.~\ref{fig:DC3D_matlab} we show the fluctuation of the potential temperature at $t=600$ s and height $z=280$ m along the line connecting the points $(0,0,280)$ and $(12000,0,280)$. Compared with the previous 2D benchmarks, the transition to a fully three-dimensional configuration results in a more challenging test case, and consequently a less pronounced agreement with the reference solution.
As in the previous tests, the solution obtained without mesh adaptation exhibits spurious oscillations, which are significantly attenuated when mesh adaptation is employed. The solution obtained with the PMA algorithm still departs noticeably from the reference profile and approaches it only in the final portion of the domain. By contrast, the IREE solution reproduces the reference trend more accurately, although visible discrepancies remain, particularly for the coarser mesh with $h=400$ m.
%

As a final comparison, we examine the computational performance of the two mesh adaptation strategies. The computational times associated with the simulations performed on meshes $h=200,400$ m are reported in Table~\ref{tab:13}, together with the corresponding percentage time savings relative to the reference simulation. 
Unlike the 2D benchmarks, where noticeable differences were observed between the computational costs of the IREE and PMA algorithms (see Tables~\ref{tab:2} and \ref{tab:10}), the gap between the two adaptive strategies becomes much less pronounced in the present 3D test. This behavior is consistent with the refinement patterns discussed above, which show a smaller discrepancy between the regions refined by the two algorithms.
At the same time, the benefits of mesh adaptation become even more significant than in the two-dimensional setting. While the finest adapted meshes in the 2D benchmarks yielded time savings of roughly 70\%, the corresponding savings in the present 3D test range from 90.7\% to 93.6\%. These figures are particularly remarkable when compared with the computational cost of the reference simulation, which requires approximately $58$ hours to complete. 
\begin{table}[h!]
\centering
\begin{tabular}{|c||c|c|}
\hline
mesh size h [m]& IREE algorithm [h]& PMA algorithm [h]  \\
\hline
200 & 5.4 (90.7\%) & 4.9 (91.5\%)\\       
400 & 4.1 (92.9\%)  &  3.7 (93.6\%)\\         
\hline
\end{tabular}
\caption{3D density current: computational performance of the IREE and PMA algorithms for mesh sizes $h = 200, 400$ m in terms of computational time and percentage of time saved with respect to the reference simulation ($58$ hours).
}\label{tab:13}
\end{table}

\section{Conclusions and future perspectives}\label{sec:conc}
We have considered an isotropic recovery-based error estimator (IREE) for mesh adaptation, a class of estimators widely used in engineering practice owing to their simplicity, robustness, and computational efficiency. 
Besides their ease of implementation and low computational overhead, recovery-based estimators possess a remarkable degree of flexibility, as they can be naturally combined with different discretization techniques. This feature has enabled us to transfer ideas and methodologies originally developed within the finite element framework to a finite volume environment.

The proposed IREE-based mesh adaptation strategy has been assessed on a set of representative atmospheric flow benchmarks. This application area is particularly relevant, since the potential of mesh adaptation for atmospheric simulations remains largely unexplored despite the maturity of adaptive methodologies in computational fluid dynamics. 

Through two classical two-dimensional atmospheric benchmarks (the rising thermal bubble and the density current) and a three-dimensional extension of the density current test case, we assessed the performance of the proposed mesh adaptation procedure against a plain mesh adaptation (PMA) algorithm which is widely adopted in the literature and implemented in several software frameworks.

Both qualitative and quantitative analyses demonstrate that mesh adaptation substantially improves the stability and accuracy of the numerical solution compared with simulations performed on non-adapted meshes, independently of the specific adaptation strategy employed. The results further indicate that the IREE algorithm is generally more effective at suppressing numerical instabilities, even when starting from relatively coarse meshes. In contrast, the PMA algorithm may still exhibit localized oscillations when relatively coarse meshes are employed. For the density current benchmarks, the IREE strategy also provides higher accuracy than PMA for the majority of the simulated time interval.

From a computational standpoint, both adaptive strategies yield significant savings with respect to simulations performed on fixed uniformly fine meshes, with reductions in computational time ranging from 35\% to 94\%. The PMA algorithm remains the less expensive option, since the IREE strategy typically refines a larger portion of the computational domain and requires the additional computation of a recovered gradient field to construct the error estimator. While this difference is clearly visible in the two-dimensional benchmarks, it becomes much less pronounced in the three-dimensional test case, where the computational costs of the two adaptive approaches are remarkably similar and both achieve time savings between 91\% and 94\%. The moderate increase in computational cost associated with IREE can therefore be regarded as the price to pay for the improved accuracy and robustness of the resulting solution.

A limitation shared by both IREE and PMA is the need to specify a number of algorithmic parameters whose values strongly influence the effectiveness of the adaptation process. In the present work, most of these parameters have been selected through a benchmark-dependent trial-and-error procedure. The automatic calibration of such parameters therefore represents an important direction for future research. Additional developments will focus on extending the proposed methodology to anisotropic mesh adaptation and on coupling spatial adaptivity with adaptive time-stepping strategies. 

\section*{Aknowledgements}
We thank Tommaso Andena, who has laid the foundations of this work during an internship at SISSA, and Yerbol Palzhanov, who has helped with the set-up of the 3D density current benchmark.\\
GR acknowledges the support provided by the European Union-NextGenerationEU, in the framework of the iNEST-Interconnected Nord-Est Innovation Ecosystem (iNEST ECS00000043 - CUP G93C22000610007) consortium. The authors would like to acknowledge also INdAM-GNCS and MIUR (Italian Ministry for University and Research) through FAREX-AROMA-CFD project, P.I. Prof. Gianluigi Rozza.
SP is member of the Gruppo Nazionale Calcolo Scientifico-Istituto Nazionale di Alta Matematica
(GNCS-INdAM) and acknowledges the INdAM–GNCS 2026 project “Metodi Numerici Avanzati per Problemi di Fluidodinamica Computazionale”. SP also
acknowledges the support
by MUR, grant Dipartimento di Eccellenza 2023–2027.

\bibliographystyle{plain}
\bibliography{bibliography}
\end{document}